\documentclass[preprint]{elsarticle}
\usepackage{hyperref}
\usepackage{afterpage}     
\usepackage{float}
\usepackage{amssymb}
\usepackage{url}
\usepackage{amsmath}
\usepackage{algorithm}
\usepackage[noend]{algpseudocode}
\usepackage{epstopdf}
\usepackage{rotating}
\usepackage{balance}
\usepackage{inputenc}
\usepackage{colortbl}
\usepackage{float}
\usepackage[caption = false]{subfig}
\usepackage{cleveref}
\usepackage{svg}
\usepackage{algorithm}
\usepackage{algorithmicx}
\usepackage{algpseudocode}
\usepackage{caption}

\biboptions{sort&compress}

\journal{}

\newcolumntype{a}{>{\columncolor[gray]{0.95}}c}

\begin{document}

\begin{frontmatter}
	
\title{Online Cluster-Based Parameter Control for Metaheuristics}

\author{Vasileios A. Tatsis}
\author{Dimos Ioannidis}

\address{Information Technologies Institute, Centre for Research and Technology Hellas (CERTH), GR-57001 Thessaloniki, Greece\\
\{vtatsis, djoannid\}@iti.gr}

	
\begin{abstract}
The concept of parameter setting is a crucial and significant process in metaheuristics since it can majorly impact their performance. It is a highly complex and challenging procedure since it requires a deep understanding of the optimization algorithm and the optimization problem at hand. In recent years, the upcoming rise of autonomous decision systems has attracted ongoing scientific interest in this direction, utilizing a considerable number of parameter-tuning methods. There are two types of methods: offline and online. Online methods usually excel in complex real-world problems, as they can offer dynamic parameter control throughout the execution of the algorithm. The present work proposes a general-purpose online parameter-tuning method called Cluster-Based Parameter Adaptation (CPA) for population-based metaheuristics. The main idea lies in the identification of promising areas within the parameter search space and in the generation of new parameters around these areas. The method's validity has been demonstrated using the differential evolution algorithm and verified in established test suites of low- and high-dimensional problems. The obtained results are statistically analyzed and compared with state-of-the-art algorithms, including advanced auto-tuning approaches. The analysis reveals the promising solid CPA's performance as well as its robustness under a variety of benchmark problems and dimensions. 
\end{abstract}
	
	
\end{frontmatter}



\section{Introduction}
\label{sec:intro}

\noindent

The explosion of artificial intelligence, along with the increasing interest in autonomous systems that can make decisions in real time, has introduced increasingly complex real-world challenges. This surge in interest has driven the scientific community toward developing automated and efficient algorithms and various parameter tuning/control methods to enhance their robustness and efficiency in such demanding applications. 

Parameter tuning/control was identified early as a critical key factor in metaheuristics as it can significantly affect their performance, convergence speed, and solution quality~\cite{EibenHM99}. Metaheuristics are based on various parameters to balance exploration and exploitation throughout the search process. Thus, selecting the appropriate parameter values can lead to substantial performance enhancements and can avoid the time-consuming and complex process of manual tuning~\cite{bergstra2012random}.

Broadly speaking, parameter tuning, known as offline tuning, focuses on finding the best static parameter values before the algorithm's execution, whereas parameter control, known as online tuning, refers to dynamically adjusting the parameter values during the algorithm's execution based on the algorithm's performance. Therefore, two main categories of parameter tuning approaches have been created: 
\begin{enumerate}[(a)]
\item \textit{Offline approaches}: The parameters are tuned prior to the application of the algorithm.
\item \textit{Online approaches}: The parameters are dynamically adapted (controlled) during the algorithm's run.
\end{enumerate}

\noindent
Offline approaches involve tuning the parameters before the algorithm is executed in a pre-processing phase driven by statistical methods. The produced outcome is reusable for the same types of problems, but they struggle to maintain adequate performance on different ones and need to reconfigure their static configurations. Various offline approaches have been proposed, including the famous Design of Experiments (DOE)~\cite{montgomery2017design}, which systematically explores the parameter space through statistical evaluation methods, and the Response Surface Methodology (RSM)~\cite{myers2016response}, which builds models to approximate the relationship between parameters and performance. Sequential Model-Based Optimization (SMBO)~\cite{hutter2011sequential} is iteratively refining a surrogate model to suggest promising parameters, while in racing algorithms, such as F-Race~\cite{birattari2002racing} and Iterated F-Race~\cite{lopez2016irace}, different parameter settings are statistically compared to promote the superior configurations and discard the inferior ones. Other notable approaches, such as CALIBRA~\cite{adenso2006fine}, combine DOE with local search, Relevance Estimation and Value Calibration (REVAC)~\cite{nannen2007efficient} estimates parameter importance to guide the search toward more promising values. ParamILS~\cite{hutter2009paramils} uses iterated local search and adaptive capping to tune the parameters, while MO-ParamILS~\cite{blot2016mo} provides multi-performance evaluation metrics to compare and improve the configurations by creating a pareto front. 

Online approaches control the parameters on the fly during the algorithm execution based on feedback from the search process. Their outcome is not reusable since it is focused on the problem at hand, but their adaptive approach to fine-tuning their behavior during the execution makes them excel in various demanding and real-world dynamic problems where static parameterization significantly fails. However, their strong dependencies with the underlying optimization algorithm obstruct the creation of general-purpose approaches, as most of them are ad-hoc procedures. Specifically, a popular generic method, such as the Adaptive Operator Selection (AOS)~\cite{fialho2010analyzing}, treats parameter values or operators as arms in a multi-armed bandit problem, updating the selection probabilities regarding the observed performance. Another generic approach is the one proposed in~\cite{joshi2020parameter}, in which parameter adaptation is based on the topological characteristics of the given optimization problem, while in~\cite{karafotias2014parameter,ting2015hybrid}, Reinforcement Learning formulates the problem of parameter control as a reinforcement learning problem, and the algorithm learns policies to map parameter adjustments with search states. 

In the context of the Differential Evolution algorithm (DE), which is also the underlying adapted algorithm of this work's proposed parameter adaptation method, many ad hoc online parameter control approaches have been proposed, focusing on adapting its parameters with significant research outcomes and distinguished performance. Currently, state-of-the-art approaches such as JADE~\cite{zhang2009jade}, jDE~\cite{TJ01}, CoDE~\cite{wang2011differential}, EPSDE~\cite{mallipeddi2011differential}, SHADE~\cite{SHADE}, and L-SHADE~\cite{L-SHADE} are adapting the parameters of DE based on various strategies and sophisticated techniques, such as encoding the parameters into individuals for coevolution or using a history of successful values that drives the parameter control. Various similar works can provide a comprehensive review and insights into parameter tuning and control in metaheuristics~\cite{Hoos11,das2016recent,tatsis2023parameter}. 

Meanwhile, recent increasing interest in artificial intelligence and its various emerging real-world applications have resulted in a considerable number of hybridized optimization algorithms with machine learning algorithms. Various techniques with sophisticated population initialization mechanisms, population reduction, and parameter setting techniques have been proposed, with promising results~\cite{sun2019survey,talbi2020machine}. For example, clustering techniques can provide more nuanced parameter control strategies since they indirectly inform parameter adjustments by distinguishing promising regions of the search space or by classifying and improving individuals or sub-populations based on cluster-specific characteristics. 

More specifically, in~\cite{dragoi2016parameter}, a cluster-based DE approach for multimodal optimization problems have been proposed, which involves clustering populations into different sub-populations to locate different optima, accompanied by a self-adaptive parameter control is employed to enhance the search ability of DE. Another cluster-based dynamic differential evolution with an external archive called CDDE$\_AR$~\cite{halder2013cluster} has been proposed, where the entire population is partitioned into several clusters according to the spatial locations of the trial solutions and evolved independently. DNDE-APC~\cite{ma2021data} proposes a novel data-driven niching DE, integrating a clustering approach, a niching technique, and a local surrogate assistant method with adaptive parameter control for the nonuniqueness of inversion. Moreover, in~\cite{bovskovic2017clustering}, a new DE for multimodal optimization has been proposed that uses self-adaptive parameter control, clustering, and crowding methods. The approach clusters the population into smaller subpopulations in such a way as to improve the algorithm's efficiency. A different approach in~\cite{kotinis2014improving} uses Mamdani-type fuzzy logic controllers (FLCs) for parameter adaptation coupled with K-medoids clustering to enable the algorithm to perform a more guided search by evolving neighboring vectors belonging to the same cluster. A similar clustering approach, namely DE-KM~\cite{kwedlo2011clustering} combining DE with K-means, has been proposed to fine-tune each candidate solution obtained by DE's mutation and crossover operations.

However, most are ad hoc procedures, and machine learning algorithms are usually used to provide key insights that facilitate subsequent parameter control steps without directly controlling them. Thus, drawing inspiration from previous works, such as the general-purpose Gradient-based Parameter Adaptation with Line Search (GPALS) method that is based on gradient approximations of the solver's performance and line search~\cite{Tatsis-GPALS}, this work also focuses on the development of a new general-purpose cluster-based online parameter tuning method for metaheuristic algorithms. The main idea of the proposed method lies in identifying promising areas in the parameter domain and iteratively generating new ones until these areas shrink and converge. More specifically, the K-means clustering algorithm was selected to identify promising areas and a custom procedure was proposed to generate the new parameters. The performance of the proposed method was verified on various benchmark problems from two established test suites, including complex low- and high-dimensional problems, while comparisons with state-of-the-art algorithms have been conducted. 

The rest of the paper is organized as follows: Section~\ref{sec:backinfo} is devoted to the background knowledge typically required to follow the proposed method that is thoroughly described in Section~\ref{sec:propmeth}. Section~\ref{sec:propDE} refers to the application of the proposed method to the differential evolution algorithm, while the experimentation assessment is presented in Section~\ref{sec:exper}, and the paper concludes in Section~\ref{sec:conclu}.

\section{Background Information}
\label{sec:backinfo}

\noindent
This section briefly describes the Differential Evolution algorithm (DE), whose parameters will be adapted through the proposed parameter adaptation method and the k-means clustering algorithm. For simplicity reasons, the optimization problem under consideration will henceforth be assumed in the general $n$-dimensional bound-constraint form.

\subsection{Differential Evolution Algorithm}

\noindent
DE is a powerful and versatile population-based optimization algorithm introduced by Storn and Price in 1997 ~\cite{storn1997differential}. Its simplicity and ability to handle various optimization problems efficiently made it among the top selections in the scientific community. DE utilizes a population $P = \left\{ x_1, x_2, \ldots, x_N \right\}$ of search agents, and each agent $x_{n}$ is a vector that represents a candidate solution within the search space $X$. The population is evolving by applying specific operations that mimic genetic procedures, such as mutation, recombination, and selection.  

\textbf{Mutation:} In DE, the mutation operation creates a mutant vector for each target vector of the population. Different mutation operators have been proposed in the bibliography. The most commonly used mutant vector variations for $\mathbf{v}_{i,g}$ for the $i$-th agent in the $g$-th generation are described as follows:

\vspace{0.2cm}
\noindent
DE/best/1:
\begin{equation}
    \mathbf{v}_{i,g} = \mathbf{x}_{best,g} + F \cdot (\mathbf{x}_{r_1,g} - \mathbf{x}_{r_2,g})
	\label{eq:mutation_1}
\end{equation}

\noindent
DE/best/2:
\begin{equation}
    \mathbf{v}_{i,g} = \mathbf{x}_{best,g} + F \cdot (\mathbf{x}_{r_1,g} - \mathbf{x}_{r_2,g}) + F \cdot (\mathbf{x}_{r_3,g} - \mathbf{x}_{r_4,g})
	\label{eq:mutation_2}
\end{equation}

\noindent
DE/rand/1:
\begin{equation}
    \mathbf{v}_{i,g} = \mathbf{x}_{r_1,g} + F \cdot (\mathbf{x}_{r_2,g} - \mathbf{x}_{r_3,g}),
	\label{eq:mutation_3}
\end{equation}

\noindent
DE/rand/2:
\begin{equation}
    \mathbf{v}_{i,g} = \mathbf{x}_{r_1,g} + F \cdot (\mathbf{x}_{r_2,g} - \mathbf{x}_{r_3,g}) + F \cdot (\mathbf{x}_{r_4,g} - \mathbf{x}_{r_5,g}),
	\label{eq:mutation_4}
\end{equation}

\noindent
DE/current-to-best/1:
\begin{equation}
    \mathbf{v}_{i,g} = \mathbf{x}_{i,g} + F \cdot (\mathbf{x}_{best,g} - \mathbf{x}_{i,g}) + F \cdot (\mathbf{x}_{r_1,g} - \mathbf{x}_{r_2,g}),
	\label{eq:mutation_5}
\end{equation}

\noindent
DE/current-to-pbest/1 (without archive):
\begin{equation}
    \mathbf{v}_{i,g} = \mathbf{x}_{i,g} + F_i \cdot (\mathbf{x}_{best_p,g} - \mathbf{x}_{i,g}) + F_i \cdot (\mathbf{x}_{r_1,g} - \mathbf{x}_{r_2,g}),
	\label{eq:mutation_6}
\end{equation}

\noindent
where $r_1$, $r_2$, $r_3$, $r_4$, and $r_5$ are randomly selected indices of vectors from the current population. The parameter $F$, defined by the user, is called the \textit{scaling factor}. It takes positive values, usually in the range $[0,1]$, and controls the amplification of the differential variation and, indirectly, the exploration abilities of the algorithm. The index $best$ denotes the best member of the population according to its function value $f$, and the index $best_p$, $p \in (0,1]$, refers to an individual randomly sampled from the best-performing $100 \, p \%$ individuals in the current population. Depending on the variation of the mutation, $F$ remains the same for all $x_i$, except for the case of Eq.~(\ref{eq:mutation_6}) where in some cases, each $x_i$ is associated with its scale factor $F_i$. 

\textbf{Crossover:} The crossover operation generates a trial vector $\mathbf{u}_{i,g}$ by mixing the components of the mutant vector $\mathbf{v}_{i,g}$ and the target vector $\mathbf{x}_{i,g}$, through binomial or exponential recombination. Specifically, the binomial recombination is defined as:
\begin{equation}
    u_{i,j,g} = 
    \begin{cases} 
      v_{i,j,g} & \text{if } \text{rand}(j) \leq CR \text{ or } j = j_{\text{rand}} \\
      x_{i,j,g} & \text{otherwise}
    \end{cases}
\end{equation}

\noindent
where $j$ indexes the components of the vector reflecting the problem's dimension $d$ and $\text{rand}(j)$ is a uniformly distributed random number function in the interval $[0,1]$. The user-defined parameter $CR$ is called the \textit{crossover rate} and controls the exploitation of the algorithm, while $j_{\text{rand}}$ is a randomly chosen index ensuring that $\mathbf{u}_{i,g}$ gets at least one component from $\mathbf{v}_{i,g}$. 

To generate the trial vector $\mathbf{u}_{i,g}$ in exponential recombination, first, the index $n$ that denotes the starting point of the crossover in the vectors involved is randomly chosen from the set of vector indices $\{1, 2, \dots, d\}$, and then the number $L$ is determined, representing the length of the crossover, which can vary between 1 and $d$. The pseudocode to obtain L is:

\begin{algorithmic}
\State $L \gets 0$
\While{$\text{rand(j)} \leq CR \text{ and } (L \leq d)$}
    \State $L \gets L + 1$
\EndWhile
\State $L \gets max\{L + 1\}$
\end{algorithmic}

\noindent
This procedure ensures that if no component is selected by chance (i.e $L=0$), then $L$ is set to 1 so that at least one component is taken from $\mathbf{v}_{i,g}$. 
Thus, for $L \geq 1$  the trial vector $\mathbf{u}_{i,g}$ is generated as follows:
\begin{equation}
    u_{i,j,g} = 
    \begin{cases} 
      v_{i,j,g} & \text{for $j \in \{((n+k -1) \mod d)+1: k=1,2, \dots, L$\}} \\
      x_{i,j,g} & \text{for all other $j \in \{1, 2, \dots, d\}.$}
    \end{cases}
\end{equation}

\textbf{Selection:} The selection process determines which vectors are carried over to the next generation by comparing the objective function values of the produced trial vector $\mathbf{u}_{i,g}$ and the target vector $\mathbf{x}_{i,g}$:
\begin{equation}
    \mathbf{x}_{i,g+1} = 
    \begin{cases} 
      \mathbf{u}_{i,g} & \text{if } f(\mathbf{u}_{i,g}) \leq f(\mathbf{x}_{i,g}) \\
      \mathbf{x}_{i,g} & \text{otherwise}
    \end{cases}
\end{equation}

\noindent
The algorithm iteratively continues until a predefined termination condition is met, such as the maximum number of function evaluations or the desired solution quality. DE's effectiveness lies firmly in its control parameters and the robust mechanisms of exploiting the search space through differential mutation combined with its straightforward selection strategy, making it suitable for various challenging optimization problems, including those with complex multimodal landscapes.

\subsection{K-Means Clustering Algorithm}
K-Means is one of the simplest and most widely used unsupervised Machine Learning (ML) algorithms, extensively used for clustering purposes. It was first proposed in 1967 by Macqueen to solve the problem of a well-known cluster~\cite{macqueen1967some}. K-Means aims at partitioning the given dataset into distinct clusters through an iterative process that tries to make the inter-cluster data points as similar as possible while keeping the clusters as far apart as possible until a stopping criterion is met. Usually, this criterion minimizes the within-cluster sum of squares (WCSS)~\cite{dao2015constrained}.

Putting it formally, let $\mathcal{D} = \{ \mathbf{d}_1, \mathbf{d}_2, \ldots, \mathbf{d}_n \}$ be a dataset in a $m$-dimensional Euclidean space $\mathbb{R}^m$, $n$ be the number of them, and $K$ be the number of clusters. 
The algorithm begins with the initialization of $K$ randomly selected cluster centers $M = (\mathbf{\mu}_1, \mathbf{\mu}_2, \ldots, \mathbf{\mu}_K)$, from the provided data set. More sophisticated methods such as K-Means++ \cite{arthur2007k} use more sophisticated initialization methods to enhance convergence. The primary goal is the minimization of the objective function: 

\[ J(M) = \sum_{j=1}^{K} \sum_{\mathbf{d}_i \in C_j} \|\mathbf{d}_i - \mathbf{\mu}_j\|^2 \]

\noindent where \(C_j\) represents the set of points in the \(j\)-th cluster and \(\mathbf{\mu}_j\) is the centroid of that cluster. Then, the algorithm follows an iterative process that consists of two main steps: 

\textbf{Assignment Step:} Each data point $\mathbf{d}_i$ is assigned to the nearest cluster determined from the Euclidean distance from the corresponding centroid. 

\[ C_i = \arg \min_{j} \|\mathbf{d}_i - \mathbf{\mu}_j\|^2 \]

\textbf{Update Step:} The centroids are updated by using the mean of all points assigned to each cluster: 

\[ \mathbf{\mu}_j = \frac{1}{|C_j|} \sum_{\mathbf{d}_i \in C_j} \mathbf{d}_i \]

\noindent These steps are iteratively executed until the centroids stabilize, indicating convergence or until a predefined number of iterations defined by the user is reached. 

K-means propose some challenges, such as the number of $K$ that needs to be defined before execution, sensitivity to the initial placement of the centroids, difficulty clustering data of varying shapes and sizes, and poor performance with high-dimensional data. Various work has been proposed to overcome these challenges, such as the Silhouette Score~\cite{rousseeuw1987silhouettes}, which measures how similar a point is to its cluster compared to the others. The silhouette score ranges from $-1$ to $+1$, and a high value indicates that the point is well matched to its group and poorly matched with the others. Similarly, many methods have also been proposed to efficiently calculate the number of clusters~\cite{pelleg2000x} and to mitigate the initialization problem~\cite{rodriguez2014clustering} by using a heuristic approach of a decision graph to identify the most promising centroids. However, despite the challenges, K-means remain a fundamental tool in unsupervised learning and data analysis due to its efficiency, simplicity, speed, and ease of implementation.


\section{Proposed Parameter Adaptation Method}
\label{sec:propmeth}

\begin{figure}[tbp]
	\centering
	\includegraphics[width=0.8\linewidth]{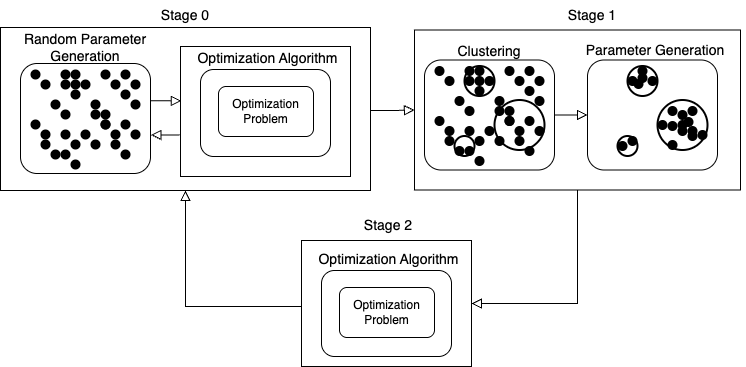}
	\caption{Graphical representation of the proposed CPA's stages.}
	\label{fig:CPA_flow}
\end{figure}

The proposed method, henceforth denoted as \textit{Cluster-based Parameter Adaptation} (CPA), draws inspiration from the nature of metaheuristic algorithms and their core mechanisms that balance the exploration and exploitation of the provided search space. The primary conception of the method was derived from the knowledge and inspiration of previous works~\cite{Tatsis2015-DEGPA, Tatsis-GPALS, tatsis2023reinforcement} that first proposed the idea of forming a new search space in the parameters domain called the \textit{parameter search space}. Consequently, since the search space of an optimization algorithm and the parameter search space share similar principles, there is plenty of space for innovative research in this direction. For this reason, extensive research with sound results has been conducted, including the application of the Particle Swarm Optimization Algorithm~\cite{Tatsis2017-PSOPNA}, the adaptation of discrete non-continuous parameters~\cite{Tatsis2016-DEGPOA}, as well as extensive sensitivity analysis~\cite{Tatsis2018-Sensitivity, Tatsis2017-IJAIT}.

The proposed method introduces three stages, namely, \textit{Random Parameter Generation}, \textit{Guided Cluster-Based Parameter Generation}, and \textit{Dynamic Deployment}. After Random Generation, the method iterates between the last two stages and stops when the computational budget is over, without spending additional resources on the parameter tuning. The proposed method proposes some extra control parameters, but as we have examined in previous works~\cite{Tatsis2017-IJAIT}, their sensitivity correlated with the overall performance of the optimization algorithms is considerably lower than directly controlling the parameters of the optimization algorithm. A graphical representation of all stages is visualized in Fig.~\ref{fig:CPA_flow}. 

\subsection{Stage 0: Random Parameter Generation}
\label{sec:randgen}
In the first stage, the underlying optimization algorithm initializes its population within the provided search space $X$ and randomly initializes its control parameters within the parameter search space $G$. Consequently, let the real-valued parameters of the optimization algorithm be denoted as:
\[ 
\xi_1,\xi_2,\ldots,\xi_{\zeta},
\]
where $\zeta$ is the number of parameters, then, the constructed $\zeta$ dimensional parameter search space is formed as follows:
\[
G = \prod_{i=1}^{\zeta} \left[ \xi_i^l, \xi_i^u \right] = \left[ \xi_1^l, \xi_1^u \right] \times \left[ \xi_2^l, \xi_2^u \right] \times \cdots \times \left[ \xi_{\zeta}^l, \xi_{\zeta}^u \right]
\]
where $l$ and $u$ represent their corresponding lower and upper bounds. At each iteration, $N$ new parameters equal to the number of population are generated randomly within $G$, using a uniform distribution:
\[
T_{j,i} \sim U(\xi_i^l, \xi_i^u), \quad \text{for } j = 1, 2, \ldots, N \quad \text{for } i = 1, 2, \ldots, \zeta.
\]
\noindent where \( U(a, b) \) denotes a uniform random variable over the interval \( [a, b] \), and in case the generated parameters \( T \) exceeds their bounds, they are set to the nearest bound:
\[
T_{j,i} = \begin{cases}
\xi_i^l, & \text{if } T_{j,i} < \xi_i^l, \\
\xi_i^u, & \text{if } T_{j,i} > \xi_i^u, \\
T_{j,i}, & \text{otherwise}.
\end{cases}
\]
The underlying optimization algorithm is executed for one iteration using the generated parameters $T$, and an external archive $A$ is maintained to store the successful parameters that produce better members in terms of fitness values compared to the corresponding old members. Initially, the archive is empty $A=\emptyset$, and the user specifies its maximum size $AS$. As the optimization progresses, successful parameters are added to $A$ until its maximum capacity. Once $A$ is full, the method proceeds to the next stage.

\subsection{Stage 1: Guided Cluster-Based Parameter Generation}
\label{sec:evalper}
This stage focuses on identifying the most promising areas of $G$ to guide parameter generation. Specifically, the K-means algorithm is applied in the archive $A$ to determine the cluster centers (centroids) that will promote the generation of new promising parameters in the \( \zeta \)-dimensional parameter space to enhance the optimization process of the underlying optimization algorithm.   

As discussed in Section~\ref{sec:backinfo}, the number of clusters $K$ in K-means should be predefined as it plays a crucial role in the algorithm's performance. Therefore, the optimal number of clusters is calculated by applying a silhouette score within the predefined range $K \in \{K_{l}, K_{u}\}$, where $l$ and $u$ represent the lower and upper limit, respectively. However, since the computation of the silhouette score is a demanding procedure, the user can set a single value for $k$. The clustering defines $k$ centroids:
\[
\mathbf{c}_k = (\xi_{k,1}, \xi_{k,2}, \ldots, \xi_{k,\zeta})
\]
where \( k = 1, 2, \ldots, K \), representing a point in the $\zeta$-dimensional parameter space. To promote a balance between exploration and exploitation, the number of new parameters to be generated for each cluster is proportionally allocated to the number of members (size) in that cluster. In this way, the more populated clusters promote exploitation, while the less populated clusters support exploration. Specifically, let $n_k$ denote the number of members in each cluster $c_k$. Then, the number of new parameters assigned to the cluster $C_k$ is determined by:
\[ 
s_k = \left\lfloor RP \times \frac{n_k}{\sum_{h=1}^{K} n_h} \right\rfloor
\]
where $RP$ denotes the total number of new parameters to be generated across all clusters. 

\subsubsection*{Random Direction Sampling}
Thus, for each centroid $c_k$, $s_k$ new parameter candidates are generated using a probability decay model that simulates the effect of evaporation as the distance from $c_k$ increases, along with a random direction sampling in the $\zeta$ dimensional space through a random unit vector $P_{k,j}$. Let $j=1,2, \dots, s_k$ index the new candidates associated with cluster $C_k$. The generation starts by sampling the unit vector:
\[
    o^{(k,j)} = \left( q_{1}^{(k,j)}, q_{2}^{(k,j)}, \ldots, q_{\zeta}^{(k,j)} \right)
\]
where each component is drawn independently from a standard normal distribution,
\[
    q_{i}^{(k,j)} \sim \mathcal{N}(0, 1) \quad \text{for } i = 1, 2, \ldots, \zeta.
\]
Then, the unit vector $P_{k,j}$ is normalized as follows:
\[
    P_{k,j} = \frac{o^{(k,j)}}{\left\| o^{(k,j)} \right\|}, \quad \text{where } \left\| o^{(k,j)} \right\| = \sqrt{\sum_{i=1}^{\zeta} \left( q_{i}^{(k,j)} \right)^2},
\]
to ensure that the directions are uniformly distributed on the surface of the $\zeta$-dimensional hyperspace and promote the unbiased exploration of $G$ in all directions. 

\subsubsection*{Distance Sampling via the Evaporation Effect}
For each new parameter, the distance \( r_{k,j} \) is determined by an evaporation-based decay model:
\[
    r_{k,j} = R \cdot \left( U_{k,j} \right)^{\alpha},
\]
where $\alpha$ is the evaporation rate, $R$ is the maximum radius, and \( U_{k,j} \sim U(0, 1) \). The evaporation rate plays a crucial role in shaping the distribution of the generated parameters within the distance from the centroid of each cluster. Specifically, if $\alpha=1$, then the generated parameters will be uniformly distributed between $[0,R]$, while when $\alpha>1$, more parameters tend to be generated closer to the centroid $c_k$, promoting exploitation, and if $0<\alpha<1$ more to the maximum radius $R$, promoting exploration. 

\subsubsection*{Offset Computation and Generation of New Parameters:}
Then, the offset from the centroid is computed as:
\[
   \Delta \xi_{k,j} = r_{k,j} \cdot P_{k,j}.
\]
The new parameter vector $p_{k,j,i}$ is obtained by adding the offset to the centroid:
\[
   p_{k,j,i} = \xi_{k,i} + \Delta \xi_{k,j,i}, \quad i = 1, 2, \ldots, \zeta.
\]
Finally, a bound check is performed:
\[
    p_{k,j,i} = \begin{cases}
    \xi_i^l, & \text{if } p_{k,j,i} < \xi_i^l, \\
    \xi_i^u, & \text{if } p_{k,j,i} > \xi_i^u, \\
    p_{k,j,i}, & \text{otherwise}.
    \end{cases}
\]
to ensure that the new parameters have been generated within the limits allowed and to discard those outside. This procedure effectively directs the generation of new parameter candidates toward the most promising regions in $G$, while balancing exploitation and exploration through random directional sampling, controlled by evaporation-based distance decay, ensuring that the new candidates are diversified and adaptively generated to the underlying distribution of promising solutions. The successful parameters (in terms of fitness values, similarly as in the previous stage) are stored in an external archive $NP$, until the archive exceeds its maximum capacity $RP$, then it goes to the next stage.

\subsection{Stage 2: Dynamic deployment}
\label{sec:dyndep}
In the final stage, the newly generated parameters $p$ are fed into the underlying optimization algorithm to guide the search process towards better solutions. Similarly to Stage 0, the underlying optimization algorithm executes one iteration at a time using the newly generated parameters until all are used. Then, the method empties the array $A$ and iterates from Stage 0 until the algorithm's maximum computational budget is reached. The corresponding pseudocode, including all the stages, is provided below.

\subsection*{CPA Pseudocode}
\label{alg:CPA}
\begin{algorithmic}[1]
\Require Number of adapted parameters $\zeta$; 
         Search space bounds SX = \(\{(\xi_j^l,\xi_j^u)\}_{j=1}^{\zeta}\); 
         Population size $N$; 
         Maximum archive size $AS$;
         Maximum archive size $RP$;
         Number of clusters  $k$;
         Radius  $R$;
         Evaporation Rate $EV$;
         Maximum number of iterations \(gMax\).
\Ensure $Solution$.

\State Initialize Population $P(SX, N)$
\State Initialize archive $A \gets \emptyset$
\State Initialize archive $NP \gets \emptyset$
\State  \(g \gets 1\)
\While{\(g \ne gMax\)}
   \If{\(A = \emptyset\)}
        \While{\(|A| \leq AS\)}
            \State \textbf{Random Parameter Generation}
            \For{\(i \gets 1\) to \(N\)}
                \For{\(j \gets 1\) to \(\zeta\)}
                    \State Sample \(T_{i,j} \sim U(\xi_j^l, \xi_j^u)\)
                    \State Apply boundary check:
                    \[
                    T_{j,i} = \begin{cases}
                    \xi_i^l, & \text{if } T_{j,i} < \xi_i^l, \\
                    \xi_i^u, & \text{if } T_{j,i} > \xi_i^u, \\
                    p_{j,i}, & \text{otherwise}.
                    \end{cases}
                    \]
                \EndFor
            \EndFor
            \State \textbf{Optimization Algorithm}
            \State \(P' \gets RunOptimization(P, T, 1)\)  // For one iteration.
            \State \(S \gets ComparePopulations(P,P')\) // Output: The successful members reflecting to the successful parameters.
            \State Append $S$ to $A$  
            \State \(P \gets P'\)
            \State \(g \gets g + 1\)
        \EndWhile
    \Else
        \State \textbf{Guided Parameter Generation}
        \For{$k \gets 1$ to $K$}
            \State Calculate the number of new parameters for cluster $C_k$:
            \[
            s_k = \left\lfloor RP \times \dfrac{n_k}{\sum_{h=1}^{K} n_h} \right\rfloor
            \]
            \For{$j \gets 1$ to $s_k$}
                \State \textbf{Random Direction Sampling} $\mathbf{P}_{k,j}$
                \For{$i \gets 1$ to $\zeta$}
                    \State Sample $q_i^{(k,j)} \sim \mathcal{N}(0, 1)$
                \EndFor
                \State Form $\mathbf{o}^{(k,j)} = \left( q_1^{(k,j)}, q_2^{(k,j)}, \ldots, q_{\zeta}^{(k,j)} \right)$
                \State Compute norm $\left\| \mathbf{o}^{(k,j)} \right\| = \sqrt{\sum_{i=1}^{\zeta} \left( q_i^{(k,j)} \right)^2}$
                \State Normalize to obtain $\mathbf{P}_{k,j} = \dfrac{\mathbf{o}^{(k,j)}}{\left\| \mathbf{o}^{(k,j)} \right\|}$
        
                \State \textbf{Distance Sampling via the Evaporation Effect} $r_{k,j}$
                \[
                r_{k,j} = R \cdot (U_{k,j})^{\alpha}, \quad \text{where } U_{k,j} \sim U(0, 1)
                \]
        
                \State \textbf{Offset Computation } $\Delta \boldsymbol{\xi}_{k,j} = r_{k,j} \cdot \mathbf{P}_{k,j}$
        
                \State \textbf{Generation of New Parameters} $\boldsymbol{p}_{k,j}$
                \For{$i \gets 1$ to $\zeta$}
                    \State $p_{k,j,i} = \xi_{k,i} + \Delta \xi_{k,j,i}$
                    \State Apply boundary check:
                    \[
                    p_{k,j,i} = \begin{cases}
                    \xi_i^l, & \text{if } p_{k,j,i} < \xi_i^l, \\
                    \xi_i^u, & \text{if } p_{k,j,i} > \xi_i^u, \\
                    p_{k,j,i}, & \text{otherwise}.
                    \end{cases}
                    \]
                \EndFor
            \EndFor
        \EndFor
        \State $NP \gets p$     // Store all new parameters to archive $NP$.
        \State $A \gets \emptyset$  // Empty the archive $A$.
        
        \State \textbf{Dynamic Deployment}
        \While{\(NP \ne \emptyset\)}
            \For{\(i \gets 1\) to \(N\)}
                \For{\(j \gets 1\) to \(\zeta\)}
                    \State \(T_{i,j} \gets NP_{i,j}\)
                \EndFor
                \State Remove \(NP_i\) from $NP$
            \EndFor
            \State \textbf{Optimization Algorithm}
            \State \(P' \gets RunOptimization(P,T, 1)\)  // For one iteration.
            \State \(P \gets P'\)
            \State \(g \gets g + 1\)
        \EndWhile
        \State \(A \gets \emptyset\)
    \EndIf
\EndWhile
\State \(Solution \gets Best(P)\)   \\ Returns the best member of the population.
\State \Return $Solution$
\end{algorithmic}

\section{Application Details on Differential Evolution}
\label{sec:propDE}

The DE algorithm was selected to demonstrate the proposed CPA method due to its parameter sensitivity and efficiency under proper tuning. The corresponding approach is henceforth denoted as CPA-DE, and the two scalar DE parameters $F$, and $CR$ are assigned in the typical range $(0,1]$, constructing the parameter space $G = (0,1] \times (0,1]$ as described in Section~\ref{sec:randgen}. Then, the method initializes the population $P$ in the search space $X$ of the corresponding underlying optimization algorithm. 

\begin{figure}[tbp]
  \centering
  \includegraphics[width=\linewidth]{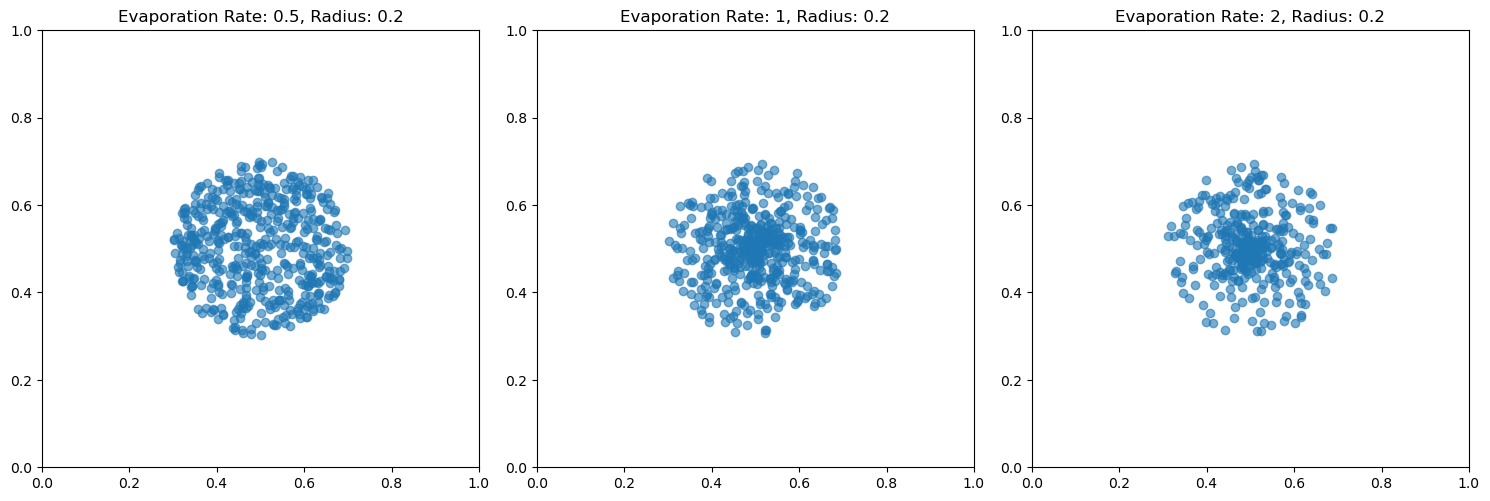}
  \caption{Graphical representation of the effect of different evaporation rates.}
  \label{fig:evap_vis}
\end{figure}

The five CPA parameters were initialized according to the default values proposed in Section~\ref{sec:propmeth} through extensive experimentation. Specifically, the maximum size of the archive $A$, the number of clusters $k$, the total number of new parameters to generate $RP$, the radius $R$, and the evaporation rate $EV$ were established at: 
\[
AS = 50, \quad k = 8 \quad RP = 200 \quad R = 0.2 \quad EV = 0.5.
\]
The most impactful parameter is the number of clusters $k$ since it can reduce the identification of promising areas in the search space $G$ that can directly impact the performance of DE. The higher the number, the more promising regions can be identified, but more processing power resources are needed. However, the number of clusters $k$ can be identified using the silhouette score or can be set manually by the user. The radius $R$ helps to obtain the convergence, as the higher the number, the more parameters are generated in a wider area around the centroid and opposite. The evaporation rate $EV$ can also play a role in the convergence of the method, as the higher the number, the generated parameters tend to be closer to the centroid. A visual representation is offered in Fig.~\ref{fig:evap_vis}. The rest of the parameters, $AS$ and $RP$, do not play a crucial role in the method's performance since they balance the resources in terms of computation terms because the less the $RP$, the more resources are needed to cluster the search space $G$. A more thorough parameter sensitivity evaluation has been conducted and the results are provided in Section~\ref{sec:exper}.

As described in Section~\ref{sec:propmeth}, the method starts by generating randomly distributed parameters within the parameter space $G$ and evaluating the population using these parameters by running the DE for one iteration. Afterwards, the successful ones are stored in the external archive $A$ until the maximum size is reached. The method then continues to stage 1, where the K-means algorithm is used to group the successful parameters of the archive $A$ and provide the centroids of each group. Through the procedure, as described in Section~\ref{sec:evalper}, $RP$ new parameters are generated and stored in the archive $NP$ to be used by the DE algorithm until the archive becomes empty, indicating the start of a new cycle. This procedure continues until the maximum function evaluations (or maximum iterations) are reached without spending any additional budget for the parameter tuning. Finally, the best-found solution was reported. 

\section{Experimental Evaluation}
\label{sec:exper}
The performance of the proposed CPA-DE approach has been validated on two established test suites that include low- and high-dimension unimodal, multimodal, hybrid, composition, and rotated and shifted test problems. The metric to measure the solution quality for all problems was the objective value error
\begin{equation}
\varepsilon^* = f\left(x^*\right) - f\left(x_{\rm opt}\right),
\label{eq:crit}
\end{equation}
\noindent where $x^*$ is the optimal known solution and $x_{\rm opt}$ the solution provided by the corresponding algorithm.

\subsection{Benchmark Problems}
\label{sec:benchprob}
Parameter tuning is essential in complex and large-scale problems where suboptimal tuning can significantly impact the algorithm's performance. Consequently, in this work, two different benchmark test suites have been selected, namely the SOCO 2011 (Soft Computing special issue on large-scale continuous optimization problems) test suite~\cite{TJ12} and the CEC 2013~\cite{CEC2013} test suite to demonstrate the efficiency and scalability of the proposed method. 

The SOCO 2011 comprises 19 continuous optimization problems scaling up to $1000$ dimensions denoted as $f_1$-$f_{19}$ with a variety of search spaces. The problems are designed to test the scalability and exploration capabilities of the tested algorithms in high-dimensional spaces that reflect some real-world scenarios. The test suite includes separable and non-separable, rotated, shifted, and hybrid functions that combine characteristics of multiple functions, including six problems, $f_1$-$f_{6}$, from the CEC 2008 test suite. The test suite specifies the maximum number of function evaluations as 
\[
\tau_{\max} = 5000 \, n,
\]
for all problems. 

Following the relevant guidelines, three base algorithms are considered the baseline for comparisons, followed by a variety of other complementary algorithms. The base algorithms include DE with exponential crossover, CHC~\cite{TJ14}, and GCMAES~\cite{TP08}, while the additional thirteen algorithms include SOUPDE~\cite{TJ15}, DE-D$^{40}$+M$^m$~\cite{TJ16}, GaDE~\cite{TJ18}, jDElscop~\cite{TJ19}, SaDE-MMTS~\cite{TJ03}, MOS~\cite{TJ21}, MA-SSW-Chains~\cite{TJ22}, RPSO-vm~\cite{TJ23}, Tuned IPSOLS~\cite{TJ24}, EvoPROpt~\cite{TJ25}, EM323~\cite{TJ26}, VXQR1~\cite{TJ27}, and GODE~\cite{TJ17}. Note that adaptive and self-adaptive algorithms are included for fair comparisons. The complementary material of the test suite is available in~\cite{TJ13}.

Along with the high-dimensional problems of SOCO 2011, the mainstream CEC 2013 test suite developed for the Special Session on Real-Parameter Single Objective Optimization at the IEEE Congress on Evolutionary Computation was also considered~\cite{CEC2013} as it is widely recognized as a comprehensive benchmark for testing real-parameter optimization algorithms. The test suite consists of $28$ highly complete optimization problems, denoted as $f_1$-$f_{28}$, including unimodal, multimodal, shifted, hybrid, and composition functions, while the search space is $X=[-100,100]^n$ for all test problems. The dimensions range from $10$ to $100$, while the maximum number of function evaluations dictated by the test suite at:
\[
\tau_{\max} = 10000 \, n.
\]
Many competitive algorithms are offered for comparisons, namely NIPOPaCMA~\cite{loshchilov2013cma} , Icmaesils~\cite{liao2013benchmark}, DRMA-LSCh-CMA~\cite{lacroix2013dynamically}, mvmo~\cite{rueda2013hybrid}, SMADE~\cite{caraffini2013super}, TLBSaDE~\cite{biswas2013teaching}, DEcfbLS~\cite{poikolainen2013differential}, b6e6rl~\cite{rueda2013hybrid}, SPSRDEMMS~\cite{zamuda2013structured}, CMAES-RIS~\cite{caraffini2013cma}, SPSOABC~\cite{el2013testing}, Jande~\cite{brest2013real}, DE APC~\cite{elsayed2013differential}, fk-PSO~\cite{nepomuceno2013self}, TPC-GA~\cite{elsayed2013genetic}, PVADE~\cite{dos2013population}, CDASA~\cite{korovsec2013continuous}, SPSO2011~\cite{zambrano2013standard}, PLES~\cite{papa2013parameter}. The complementary material from the test suite is also available online at~\cite{CEC_2013_Material_Url,SUG_Url}, while it shall be noted that the competitor algorithms adopt the suggested parameter setting as proposed in their original sources, and the provided computational budget was not reduced due to their \textit{a priori} tuning.

\subsection{Parameter Sensitivity Results for the SOCO 2011 test suite}
\label{sec:par_sens_SOCO}

\begin{table}[t]
	\footnotesize
	\renewcommand{\tabcolsep}{8pt}
	\renewcommand{\arraystretch}{1.5}
	\begin{center}
		\caption{Statistical comparisons of various CPA-DE variants and the baseline CPA-DE variant on the SOCO 2011 test suite.}
		\label{tab:compGPALS_SOCO}
		\begin{tabular}{ccccc}
			\hline
			& & \multicolumn{3}{c}{CPA$_{4\_0.05\_100\_100}-DE_{R}^{60}$} \\
			\cline{3-5} 
			Dim. & Algorithm & $+$ & $-$ & $=$ \\
			\hline
			$50$ 
			& CPA$_{2\_0.05\_100\_100}-DE_{R}^{60}$      & $2$ & $4$ & $13$ \\
			& CPA$_{6\_0.05\_100\_100}-DE_{R}^{60}$      & $3$ & $0$ & $16$ \\
			& CPA$_{8\_0.05\_100\_100}-DE_{R}^{60}$      & $5$ & $0$ & $14$ \\
			\hline
   
			& CPA$_{4\_0.1\_100\_100}-DE_{R}^{60}$       & $1$ & $0$ & $18$ \\
			& CPA$_{4\_0.15\_100\_100}-DE_{R}^{60}$ 	 & $3$ & $2$ & $14$ \\
			& CPA$_{4\_0.2\_100\_100}-DE_{R}^{60}$       & $5$ & $2$ & $12$ \\
			\hline
   
            & CPA$_{4\_0.05\_50\_100}-DE_{R}^{60}$       & $4$ & $1$ & $14$ \\
			& CPA$_{4\_0.05\_150\_100}-DE_{R}^{60}$ 	 & $1$ & $2$ & $16$ \\
			& CPA$_{4\_0.05\_200\_100}-DE_{R}^{60}$ 	 & $0$ & $3$ & $16$ \\
			\hline
   
            & CPA$_{4\_0.05\_100\_50}-DE_{R}^{60}$       & $2$ & $3$ & $14$ \\
			& CPA$_{4\_0.05\_100\_150}-DE_{R}^{60}$      & $2$ & $1$ & $16$ \\
			& CPA$_{4\_0.05\_100\_200}-DE_{R}^{60}$      & $5$ & $0$ & $14$ \\
			\hline
		\end{tabular}
	\end{center}
\end{table}

The proposed CPA method introduces $5$ additional control parameters that can affect its performance. However, as has been shown from previous work~\cite{Tatsis2017-IJAIT, Tatsis2018-Sensitivity}, the additional parameters introduced by advanced tuning techniques like the proposed one, often exhibit a reduced sensitivity compared to the parameters of the original algorithm, resulting in a reduced impact on overall performance. These findings have also been validated in other works like~\cite{ojha2022assessing}, where the authors demonstrated that only certain hyperparameters have more pronounced effects on performance, suggesting that the rest can be configured more easily and with less precision without significantly affecting results. Therefore, identifying the most influential parameters of the tuning process can lead to more efficient results, as less sensitive parameters can exert less impact on the final outcome. Thus, a preliminary sensitivity analysis was conducted by individually analyzing the impact of each parameter on the overall performance of the method. 

In this work, the same settings as those proposed in previous work ~\cite{Tatsis2015-DEGPA, Tatsis-GPALS}, regarding the underlying optimization algorithms and specifically the DE algorithm, have also been adopted. Specifically, the DE/rand/1 mutation operator of Eq.~(\ref{eq:mutation_3}) was used along with exponential recombination and a fixed population size $N=60$. Each of the $4$ additional control parameters except the parameter $EV$ was discretized, and a reasonable set was initialized according to our preliminary knowledge and specifically to
\[
AS = \{50, 100, 150, 200\}, \quad k = \{2, 4, 6, 8\},
\]
\[
\quad RP = \{50, 100, 150, 200\}, \quad R = \{0.05, 0.1, 0.15, 0.2\}
\]
for all the parameters. We decided to exclude the parameter $EV$ from the sensitivity analysis, since this parameter is responsible for spreading the randomly generated points more uniformly within the radius $R$, as we can observe in Fig.~\ref{fig:evap_vis}. Then, alternating one parameter at a time, as visualized in Table~\ref{tab:compGPALS_SOCO}, the produced CPA-DE variants denoted as $CPA_{k\_R\_AS\_RP}-DE_{R}^{60}$ were compared to a baseline CPA variant. The baseline $CPA_{4\_0.05\_100\_100}-DE_{R}^{60}$ was specified based on two key reasons: empirical performance superiority as well as logical parameter selection. Specifically, this variant showed superior performance in our preliminary experiments in all tested problem instances, establishing it as a reliable benchmark for further performance assessment. Furthermore, its parameters were selected based on their logical coherence with the algorithm's theoretical framework, its design internal characteristics, and the underlying core mechanisms in an attempt to balance performance and computation resources. 

As we observe in Table~\ref{tab:compGPALS_SOCO}, the analysis was performed in the lowest dimension of the SOCO 2011 test suite due to the low computational execution time. The first three variants refer to the number of clusters $k$. The best-performing variant is the one with the highest number of clusters and specifically with $k=8$, which is also in line with the logical explanation that the higher the number of clusters, the better the exploration of the parameter space and the generation of new promising parameters accordingly. However, since the rise of clusters can affect the execution time of the algorithm for this work, we decided to keep it within these predefined bounds. The following parameter is the radius $R$; as we can see, higher values are also preferred. This selection can also be logically explained since the higher the radius of the selected clusters, the better the exploration of this region, which is a clear advantage in the early stages of almost any optimization algorithm []. 

The third parameter, $AS$, refers to the maximum size of the archive $A$, following the opposite direction to the lowest bound and specifically at the value $50$. This tendency probably occurred because parameter tuning benefits more from more frequent cluster updates that can increase the possibility of improving the population. The last parameter, $RP$, which is the total number of newly generated parameters, appears to also affect the performance of the method when higher values are selected, specifically at $200$. Thus, summarizing the findings of this preliminary sensitivity analysis, we concluded that the most effective variant of CPA (compared to baseline) is $CPA_{8\_0.2\_50\_200}-DE_{R}^{60}$ and thus adopted it for the selected benchmark problems.

\subsection{Results for the SOCO 2011 test suite}
\label{sec:SOCO}

\begin{table}[t]
	\footnotesize
	\renewcommand{\tabcolsep}{8pt}
	\renewcommand{\arraystretch}{1.2}
	\begin{center}
		\caption{Statistical comparisons between CPA-DE and the base algorithms on the high-dimensional SOCO 2011 test suite.}
		\label{tab:compBaseSOCO}
		\begin{tabular}{clrrrrrrr}
			\hline
			& & \multicolumn{3}{c}{CPA-DE$_{R}^{60}$} & & \multicolumn{3}{c}{CPA-DE$_{C}^{60}$} \\
			\cline{3-5} \cline{7-9} 
			Dim. & Algorithm & $+$ & $-$ & $=$ & & $+$ & $-$ & $=$  \\
			\hline
			$50$ & DE$_{\rm bin}$ 	& $14$ & $3$ & $2$ & & $7$  & $11$ & $1$ \\
			& DE$_{\rm exp}$ 		& $11$ & $5$ & $3$ & & $4$  & $13$ & $2$ \\
			& CHC 	       			& $19$ & $0$ & $0$ & & $16$ & $1$  & $2$ \\
			& GCMAES 				& $16$ & $3$ & $0$ & & $15$ & $3$  & $1$ \\
			\hline
			$100$ & DE$_{\rm bin}$ 	& $17$ & $1$ & $1$ & & $7$  & $10$ & $2$ \\
			& DE$_{\rm exp}$ 		& $13$ & $6$ & $0$ & & $3$  & $15$ & $1$ \\
			& CHC 	       			& $19$ & $0$ & $0$ & & $13$ & $3$  & $3$ \\
			& GCMAES 				& $16$ & $3$ & $0$ & & $15$ & $4$  & $0$ \\
			\hline
			$200$ & DE$_{\rm bin}$ 	& $17$ & $1$ & $1$ & & $7$  & $10$ & $2$ \\
			& DE$_{\rm exp}$ 		& $13$ & $6$ & $0$ & & $3$  & $16$ & $0$ \\
			& CHC 	       			& $19$ & $0$ & $0$ & & $12$ & $2$  & $5$ \\
			& GCMAES 				& $16$ & $3$ & $0$ & & $14$ & $4$  & $1$ \\
			\hline
			$500$ & DE$_{\rm bin}$ 	& $19$ & $0$ & $0$ & & $11$ & $6$ & $2$ \\
			& DE$_{\rm exp}$ 		& $12$ & $7$ & $0$ & & $1$  & $17$ & $1$ \\
			& CHC 	       			& $19$ & $0$ & $0$ & & $9$  & $8$  & $2$ \\
			& GCMAES 				& n/a  & n/a & n/a & & n/a  & n/a  & n/a \\
			\hline
		\end{tabular}
	\end{center}
\end{table}

The CPA-DE variant with the best performance as examined through the preliminary sensitivity analysis of Section~\ref{sec:par_sens_SOCO} adopted for experimentation. In addition to the examined DE mutation operator used for the sensitivity analysis, the DE/current-to-pbest/1 with external archive was also considered since this operator has shown significant performance improvements in various problems. Consequently, the subscripts $R$ and $C$ denote the type of mutation operator for Rand/1 and current-to-pbest/1 with the external archive accordingly, while the superscript denotes the population number. 

Following the test suite guidelines, the performance assessment was carried out in three stages, while $25$ independent experiments were conducted for all the competitors. In the first stage, CPA-DE was compared with three baseline algorithms proposed by the test suite: DE with exponential/binomial recombination and the Rand/1 mutation operator, CHC, and GCMAES, using the optimal parameters proposed in their original articles. Comparisons of the second stage were made against various optimization algorithms and precisely $13$ DE- and non-DE-based algorithms, employing various sophisticated procedures and self-adapted mechanisms. The final stage of the comparisons includes the state-of-the-art SHADE~\cite{SHADE} algorithm, which is constituted as a key pillar of inspiration in the creation of more sophisticated improvements with superior results in recent years. The reason behind the selection of the SHADE algorithm and not of the newer versions was that the adaptation mechanisms that SHADE incorporates remained the same in the next improved versions, and the comparison with the proposed method would be unfair if a newer variant were selected since we both focus explicitly on the adaptation of the DE's scalar parameters.

The results of the Wilcoxon Ranksum tests at the confidence level $95\%$ of the first phase against the base algorithms are presented in Table~\ref{tab:compBaseSOCO}. The total number of wins, losses, and draws for a pair of algorithms provided from the statistical tests are denoted by the symbols ``$+$'', ``$-$'', and ``$=$'', respectively. These results provide a comparative statistical analysis of the CPA-DE algorithm compared to the baseline algorithms of the test suite in various dimensions up to $500$. Two variants of CPA-DE: CPA-DE$_{R}^{60}$ and CPA-DE$_{C}^{60}$ were selected to provide an in-depth overview of the performance differences between the two selected mutation operators. As we can observe, both variants show competitive performance, winning a significant number of test cases compared to the baseline algorithms. In particular, CHC consistently performs the worst, with zero wins compared to the DE$_{R}^{60}$ variant, while DE$_{C}^{60}$ achieves slightly better results. For the highest dimension, the results for GCMAES are marked as "n/a", indicating that the evaluation was not conducted within reasonable time limits for this instance. Thus, the consistent performance of CPA-DE$_{R}^{60}$ highlights its robustness across different problem dimensions, particularly in high-dimensional optimization problems, making it a promising choice for complex optimization problems. 

The improved performance of CPA-DE$_{R}^{60}$ compared to CPA-DE$_{C}^{60}$ can be attributed to several key factors related to fundamental differences in the mutation operators. Specifically, the DE/rand/1 mutation operator is known for its strong exploration capabilities, as it randomly selects individuals from the population to generate new candidate solutions. This randomness provides an edge regarding the exploration of the search space and also avoids premature convergence to local optima, which can be particularly beneficial in high-dimensional problems where the landscape consists of many local optima. In contrast, the DE/current-to-pbest/1 with external archive guides the search towards promising regions based on the best individuals found so far (archived individuals), providing a more exploitable approach. Thus, this tendency can potentially lead to a loss of diversity, promote stagnation in high-dimensional complex search spaces, and be more prone to noise since the best historical solutions stored in the archive can mislead the search. The average errors and standard deviations obtained from CPA-DE and the base algorithms are provided in Tables~\ref{tab:SOCO_50}-\ref{tab:SOCO_500} in~\ref{sec:app1}. 

In the second phase of comparisons, the CPA-DE$_{R}^{60}$ variant was compared against $13$ additional competitors, as Table~\ref{tab:compOtherSOCO} indicates. The same notation was followed here for wins, losses, and draws, but the comparisons were conducted with the obtained average errors. CPA-DE$_{R}^{60}$ variant performs comparably with DE-based sophisticated algorithms, such as SOUPDE, GODE, GaDE, jDElscop, and SaDE-MMTS, in lower dimensions and slightly outperforms them in higher ones. This suggests that the competent algorithms can also achieve balanced exploration and exploitation in most cases except in high dimensions where CPA-DE has an advantage. The only algorithm that outperforms CPA-DE is the MOS algorithm, which is a multiple offspring sampling-based dynamic memetic DE that combines the explorative/exploitative strength of two heuristic search methods that separately obtain very competitive results in either low- or high-dimensional search spaces. For the rest of the algorithms, such as MA-SSW-Chains, EM323, and VXQR1, the same effect of increasing performance is observed as the dimension increases for CPA. Furthermore, it is worth mentioning that all the competitive algorithms are extensively fine-tuned for the specific test suite and used with the parameters as proposed in their original papers. The corresponding average errors of the competing algorithms are reported in Tables~\ref{tab:SOCO11_Errors1} and~\ref{tab:SOCO11_Errors4} in~\ref{sec:app1}. 

The last phase of the comparisons consists of head-to-head comparisons between CPA-DA and the state-of-the-art SHADE algorithm. Since in the original paper, SHADE incorporates the DE/current-to-pbest/1 with an external archive mutation operator along with its proposed adaptation mechanisms without using the mutation operator DE/rand/1, it would be unfair to compare its performance without the mutation operator dictated by the test suite and also used in CPA-DE. For this reason, all possible combinations stemming from the different population sizes and the mutation operators were considered for an in-depth comparison as visualized in Table~\ref{tab:RLGPALS_SHADE_SOCO}. The same notation with the subscript to denote the mutation operator and the superscript to denote the population size was adopted here, and the Wilcoxon Ranksum tests at the confidence level $95\%$ were used to evaluate the performance differences. Thus, the CPA-DE variants using the DE/rand/1 (denoted with ``$R$'') exhibited superior performance for both population sizes $60$ and $100$ against the corresponding SHADE approaches even with the variants using the DE/current-to-pbest/1 with external archive (denoted with ``$C$''). The only variants in which SHADE achieved better performance were the variants with the DE/current-to-pbest/1 with external archive, proving that the adaptation mechanisms of SHADE probably favor this operator. However, when the proposed $CPA-DE_{R}$ compared with the original $SHADE_{C}$ ( DE/current-to-pbest/1 with external archive), CPA-DE outperformed significantly SHADE, indicating that the adaptation mechanisms of CPA-DE are superior compared to the ones proposed by SHADE in the complex high-dimensional search spaces such as the ones proposed in this test suite, where an optimal balance of exploration and exploitation is required to maintain solid performance. The average errors for each approach are reported in Tables~\ref{tab:SHADE-SOCO-ERRORS1}-~\ref{tab:SHADE-SOCO-ERRORS2} in~\ref{sec:app1}.

\noindent


\begin{table}[t]
	\footnotesize
	\renewcommand{\tabcolsep}{9pt}
	\renewcommand{\arraystretch}{1.2}
	\begin{center}
		\caption{Comparisons of the average solutions between CPA-DE$_{R}^{60}$ and a variety of algorithms on the high-dimensional SOCO 2011 test suite.}
		\label{tab:compOtherSOCO} 
		\begin{tabular}{lcrrrcrrr}
			\hline 
			& & \multicolumn{7}{c}{Dimension} \\
			\cline{3-9}
			& & \multicolumn{3}{c}{50} & & \multicolumn{3}{c}{100} \\
			\cline{3-5} \cline{7-9} 
			& & $+$ & $-$ & $=$ & & $+$ & $-$ & $=$	\\
			\hline
			SOUPDE             & &    5 & 5 & 9     & &   6 & 6 & 7 \\
			DE-D$^{40}$+M$^m$  & &   10 & 4 & 5     & &  10 & 5 & 4 \\
			GODE               & &    7 & 6 & 6     & &   7 & 6 & 6 \\
			GaDE               & &    6 & 5 & 8     & &   4 & 7 & 8 \\
			jDElscop           & &    1 & 7 & 11    & &   3 & 8 & 8 \\
			SaDE-MMTS          & &    4 & 6 & 9     & &   3 & 7 & 9 \\
			MOS 			   & &    1 & 5 & 13    & &   0 & 8 & 11 \\
			MA-SSW-Chains      & &    15 & 3 & 1    & &  15 & 4 & 0 \\
			RPSO-vm 		   & &    14 & 2 & 3    & &   13 & 2 & 4 \\	
			Tuned IPSOLS 	   & &    11 & 5 & 3    & &   10 & 6 & 3 \\
			EvoPROpt		   & &    18 & 1 & 0    & &   18 & 1 & 0 \\	
			EM323 		       & &    11 & 3 & 5    & &   11 & 4 & 4 \\		
			VXQR1 			   & &    10 & 5 & 4    & &   12 & 4 & 3 \\
			\hline
			& & \multicolumn{3}{c}{200} & & \multicolumn{3}{c}{500} \\
			\cline{3-5} \cline{7-9} 
			& & $+$ & $-$ & $=$ & & $+$ & $-$ & $=$	\\
			\hline
			SOUPDE             & &    9 & 5 & 5    & &  9 & 5 & 5 \\
			DE-D$^{40}$+M$^m$  & &    9 & 5 & 5    & &  10 & 5 & 4 \\
			GODE               & &    8 & 5 & 6    & &  8 & 5 & 6 \\
			GaDE               & &    7 & 6 & 6    & &   7 & 6 & 6 \\
			jDElscop           & &    5 & 6 & 8     & &  5 & 6 & 8 \\
			SaDE-MMTS          & &    5 & 5 & 9    & &   5 & 5 & 9 \\
			MOS 			   & &    0 & 7 & 11   & &   0 & 7 & 12 \\
			MA-SSW-Chains      & &    18 & 1 & 0   & &   18 & 1 & 0 \\
			RPSO-vm 		   & &    14 & 2 & 3   & &   14 & 3 & 2 \\	
			Tuned IPSOLS 	   & &    9 & 6 & 4    & &   9 & 6 & 4 \\
			EvoPROpt		   & &    17 & 1 & 1   & &   17 & 1 & 1 \\	
			EM323 		       & &    13 & 2 & 4   & &   13 & 2 & 4 \\		
			VXQR1 			   & &    12 & 4 & 3   & &   12 & 4 & 3 \\
			\hline
		\end{tabular}
	\end{center}
\end{table}

\afterpage{\clearpage}

\begin{table}[t]
	\footnotesize
	\renewcommand{\tabcolsep}{9pt}
	\renewcommand{\arraystretch}{1.2}
	\begin{center}
		\caption{Statistical comparisons between CPA-DE and the state-of-the-art SHADE algorithm on the SOCO 2011 test suite.}
		\label{tab:RLGPALS_SHADE_SOCO} 
		\begin{tabular}{lcrrrcrrr}
			\hline 
			& & \multicolumn{7}{c}{Dimension} \\
			\cline{3-9}
			& & \multicolumn{3}{c}{50} & & \multicolumn{3}{c}{100} \\
			\cline{3-5} \cline{7-9} 
			& & $+$ & $-$ & $=$ & & $+$ & $-$ & $=$	\\
			\hline
			CPA-DE$_{R}^{60}$  vs SHADE$_{R}^{60}$  & & 12 & 6 & 1 & & 12 & 5 & 2 \\
			CPA-DE$_{C}^{60}$  vs SHADE$_{C}^{60}$ & & 4 & 12 & 3 & & 2 & 13 & 4 \\
			CPA-DE$_{R}^{60}$  vs SHADE$_{C}^{60}$ & & 11 & 6 & 2 & & 15 & 2 & 2 \\
			CPA-DE$_{C}^{60}$  vs SHADE$_{R}^{60}$ & & 4 & 14 & 1 & & 3 & 16 & 0 \\
			CPA-DE$_{R}^{100}$  vs SHADE$_{R}^{100}$  & & 17 & 1 & 1 & & 15 & 4 & 0 \\
			CPA-DE$_{C}^{100}$  vs SHADE$_{C}^{100}$ & & 6 & 10 & 3 & & 2 & 14 & 3 \\
			CPA-DE$_{R}^{100}$  vs SHADE$_{C}^{100}$ & & 4 & 14 & 1 & & 13 & 6 & 0 \\
			CPA-DE$_{C}^{100}$  vs SHADE$_{R}^{100}$ & & 12 & 6 & 1 & & 3 & 15 & 1 \\
			\hline
			& & \multicolumn{3}{c}{200} & & \multicolumn{3}{c}{500} \\
			\cline{3-5} \cline{7-9} 
			& & $+$ & $-$ & $=$ & & $+$ & $-$ & $=$	\\
			\hline
			CPA-DE$_{R}^{60}$  vs SHADE$_{R}^{60}$  & & 14 & 5 & 0 & & 18 & 1 & 0 \\
			CPA-DE$_{C}^{60}$  vs SHADE$_{C}^{60}$ & & 4 & 12 & 3 & & 3 & 14 & 2 \\
			CPA-DE$_{R}^{60}$  vs SHADE$_{C}^{60}$ & & 18 & 1 & 0 & & 18 & 1 & 0 \\
			CPA-DE$_{C}^{60}$  vs SHADE$_{R}^{60}$ & & 3 & 16 & 0 & & 1 & 18 & 0 \\
			CPA-DE$_{R}^{100}$  vs SHADE$_{R}^{100}$  & & 15 & 4 & 0 & & 18 & 1 & 0 \\
			CPA-DE$_{C}^{100}$  vs SHADE$_{C}^{100}$ & & 2 & 14 & 3 & & 1 & 14 & 4 \\
			CPA-DE$_{R}^{100}$  vs SHADE$_{C}^{100}$ & & 13 & 6 & 0 & & 18 & 1 & 0 \\
			CPA-DE$_{C}^{100}$  vs SHADE$_{R}^{100}$ & & 3 & 15 & 1 & & 1 & 18 & 0 \\
			\hline
		\end{tabular}
	\end{center}
\end{table}

\afterpage{\clearpage}


\subsection{Results for the CEC 2013 test suite}
\label{sec:CEC}

The performance evaluation of CPA-DE for the CEC 2013 test suite followed the same experimental design as in the previous test suite. In the first phase, various optimization algorithms provided by the test suite were used for the comparisons, followed by the second phase, in which a head-to-head comparison with SHADE was conducted. All tests were based again on the Wilcoxon rank sum test at the confidence level $95\%$. 

Starting with the results of the first phase, the best performing variant CPA-DE$_{R}^{60}$ was also adopted for the comparisons without modifying any of its settings to provide a comprehensive overview of the method's adaptability to various problems without further fine-tuning. In contrast, as proposed in their original papers for this specific test suite, the competent algorithms were used again with their fine-tuned parameters. Moreover, since the specific DE/current-to-pbest/1 mutation operator with archive proved to be really effective for the problems of this test suite, it was also adopted in the CPA-DE experimentation. 

The results of Table~\ref{tab:CEC13} indicate that CPA-DE$_{R}^{60}$ can achieve competitive and, in some cases, superior results even compared to various sophisticated and complicated optimization algorithms. For example, compared to DE-based variants such as SMADE, TLBSaDE, SPSRDEMMS, and SDEcfbLS, cPA-De achieved competitive results, while with jande, DE\_APC, and PVADE the performance differences were more prominent for both dimensions, indicating that the tuning mechanisms behind CPA-DE can provide robust adaptation even when compared with algorithms that incorporate multiple mutation strategies, increasing and decreasing population sizes, local search, and ad hoc self-adaptive schemes. Compared with non-DE-based algorithms such as b6e6rl, CMAES-RIS, SPSOABC, fk-PSO, TPC-GA, CDASA, SPSO2011 and PLES, employing various types of crossover, competitive adaptive schemes, fitness-based metrics, and hybridization techniques, the results were superior, achieving in some cases wins for more than $85\%$ of the problems. The only algorithms that CPA tend to struggle with are algorithms that incorporate sophisticated local search techniques (e.g icmaesils), and novel population-based meta-heuristic techniques employing special mapping functions based on mean and variance metrics (e.g  mvmo), but even in these cases, the results are promising and paves the way for future improvements. The strong presence of CPA-DE in a completely different test suite without further modification proves the robustness and adaptability of the internal mechanisms to various optimization problems without burdening the end user. The average errors for each approach are reported in Tables~\ref{tab:CEC13-ERRORS1}-~\ref{tab:CEC13-ERRORS3} in~\ref{sec:app2}.

Moving to the second phase of comparisons with the state-of-the-art SHADE, the same experimentation methodology used in the previous test suite was also used here for fair comparisons, and the results are reported in Table~\ref{tab:CPA-DE_SHADE_CEC13}. As is clearly observable, CPA-DE$_{R}$ excels over SHADE$_{R}$ for all population sizes and dimensions, while the opposite is observed with CPA-DE$_{C}$ and SHADE$_{C}$. When we compared the two best-performing variants, namely CPA-DE$_{R}$ and SHADE$_{C}$, CPA-DE struggled to maintain solid performance within the specified budget dictated by the test suite. These performance differences verify our previous findings that the adaptation mechanisms of SHADE probably favor the DE/current-to-pbest/1 with archive mutation operator and the problems of the specific test suite. For example, it is worth examining the fact that the search spaces of this test suite lie in the range of $[-100, 100]$, while SOCO's search spaces vary, and this may have a significant impact on the robustness of the internal adaptation mechanisms to changes in domain scale. The average errors for each approach are reported in Tables~\ref{tab:SHADE-CEC13-ERRORS1}-~\ref{tab:SHADE-CEC13-ERRORS3} in~\ref{sec:app2}.

\begin{table}[t]
	\footnotesize
	\renewcommand{\tabcolsep}{10pt}
	\renewcommand{\arraystretch}{1.2}
	\begin{center}
		\caption{Statistical comparisons among CPA-DE and a variety of algorithms on the CEC 2013 test suite.}
		\label{tab:CEC13} 
		\begin{tabular}{lcrrrcrrr}
			\hline 
			& & \multicolumn{7}{c}{Dimension} \\
			& & \multicolumn{3}{c}{30} & & \multicolumn{3}{c}{50} \\
			\cline{3-5} \cline{7-9} 
			CPA-DE$_{R}^{60}$ vs. & & $+$ & $-$ & $=$ & & $+$ & $-$ & $=$	\\
			\hline
                icmaesils       & &		9    & 17    &  2   & &	    8   &  16   &  4	\\
                DRMA-LSCh-CMA   & &		6    & 17    &  5   & &		8   &  17   &  3	\\
                NIPOPaCMA       & &		10   & 14    &  4   & &		10   & 12   &  6	\\
			mvmo		 	& &		6    & 18    &  4   & &		6   &  19   &  3	\\
			SMADE	 		& &     11   & 10    &  7	& &		12   &  9   &  7	\\
			TLBSaDE			& & 	10   & 12    &  6	& &		10  &  10    & 8	\\
			DEcfbLS   		& & 	10   &  9   &  9	& &		12  &  10   &  6	\\
			b6e6rl	 		& & 	15   &  7   &  9	& &		12  &  7    &  9	\\
			SPSRDEMMS	  	& & 	12   &  8   &  8	& &		10  &  8   &  10	\\
			CMAES-RIS		& & 	12   &  8   &  8	& &		13  &  10   &  5	\\
			SPSOABC	 		& & 	15   &  6   &  7	& &		15  &  7    &  6	\\
			jande		 	& & 	12   &  8  &  8	    & &		14  &   9   &  5	\\
			DE\_APC			& & 	14   &  9   &  5	& &		14  &   9   &  5	\\
			fk-PSO	 		& & 	16   &  6   &  6	& &		10  &   8   &  7	\\
			TPC-GA	 		& & 	17   &  6   &  5	& &		16  &   8   &  4	\\
			PVADE			& & 	15   &  9   &  4	& &		18  &   8   &  2	\\
			CDASA			& & 	21   &  4   &  3	& &		21  &   5   &  2	\\
			SPSO2011		& & 	21   &  3   &  4	& &		22  &   3   &  3	\\
			PLES			& & 	24   &  1   &  3	& &		24  &   1   &  3	\\
			\hline
		\end{tabular}
	\end{center}
\end{table}

\begin{table}[t]
	\footnotesize
	\renewcommand{\tabcolsep}{9pt}
	\renewcommand{\arraystretch}{1.2}
	\begin{center}
		\caption{Statistical comparisons between CPA-DE and the state-of-the-art SHADE algorithm on the CEC 2013 test suite.}
		\label{tab:CPA-DE_SHADE_CEC13} 
		\begin{tabular}{lcrrrcrrr}
			\hline 
			& & \multicolumn{7}{c}{Dimension} \\
			\cline{3-9}
			& & \multicolumn{3}{c}{30} & & \multicolumn{3}{c}{50} \\
			\cline{3-5} \cline{7-9} 
			& & $+$ & $-$ & $=$ & & $+$ & $-$ & $=$	\\
			\hline
			CPA-DE$_{R}^{60}$  vs SHADE$_{R}^{60}$  & & 15 & 7 & 6 & & 15 & 8 & 5 \\
			CPA-DE$_{C}^{60}$  vs SHADE$_{C}^{60}$ & & 4 & 18 & 6 & & 4 & 21 & 3 \\
			CPA-DE$_{R}^{60}$  vs SHADE$_{C}^{60}$ & & 8 & 16 & 4 & & 6 & 19 & 3 \\
			CPA-DE$_{C}^{60}$  vs SHADE$_{R}^{60}$ & & 14 & 9 & 5 & & 8 & 17 & 3 \\
			CPA-DE$_{R}^{100}$  vs SHADE$_{R}^{100}$  & & 12 & 8 & 8 & & 11 & 11 & 6 \\
			CPA-DE$_{C}^{100}$  vs SHADE$_{C}^{100}$ & & 6 & 17 & 5 & & 9 & 17 & 2 \\
			CPA-DE$_{R}^{100}$  vs SHADE$_{C}^{100}$ & & 8 & 18 & 2 & & 8& 18 & 2 \\
			CPA-DE$_{C}^{100}$  vs SHADE$_{R}^{100}$ & & 12 & 10 & 6 & & 11 & 11 & 6 \\
			\hline
		\end{tabular}
	\end{center}
\end{table}

\noindent


\section{Conclusions}
\label{sec:conclu}

\noindent
Proper parameter setting is a crucial procedure that must be performed carefully and with the proper expertise in metaheuristic algorithms, since it can significantly affect their overall performance and effectiveness. Proper parameter tuning directly influences the balance between exploration and exploitation, convergence speed, and solution quality. Parameter tuning usually requires a considerable amount of time and effort, and in time-critical real-world applications, this can be a significant drawback. 

Therefore, this paper tackles the parameter tuning problem by proposing an online method that can dynamically adapt the parameters on the fly through the application of machine learning clustering methods such as K-Means to identify the most promising parameters to guide the search towards them. The proposed method achieved a significant performance leap with a variety of optimization algorithms, even with some state-of-the-art competitors that have revolutionized the area of hybrid parameter adaptation algorithms, such as SHADE. 

Both our proposed CPA and SHADE's state-of-the-art parameter adaptation schemes exhibit on-par performance in the benchmarked test problems used in this work. However, our method achieved slightly better overall performance than SHADE since, in the SOCO 2011 high-dimensional test suite, it fell behind CPA even when compared to the best SHADE variant. This observation verifies some of our previous work and sheds light on future work on the intriguing behavior of these two competitors under various circumstances and problems. Finally, the results of this work increased our motivation for future investigation of hybridizing machine learning methods with optimization-based mechanisms to further improve the proposed parameter control method and also create new ones for optimization algorithms.


\section*{Acknowledgment}

This research was funded by the European Commission through the project HORIZON EUROPE - INNOVATION ACTIONS (IA) - grant number 101058453 - FLEXIndustries. The opinions expressed in this paper are those of the authors and do not necessarily reflect the views of the European Commission.


\bibliographystyle{elsarticle-num}
\bibliography{bibliography_tatsis,bib_general}


\newpage
\appendix

\section{Results for SOCO 2011}
\label{sec:app1}

\begin{table}[t]
	\scriptsize
	\renewcommand{\tabcolsep}{5.0pt}
	\renewcommand{\arraystretch}{1.5}
	\begin{center}
		\caption{Average errors and standard deviations of the CPA-DE$_{R}^{60}$ and the SOCO base algorithms for dimension $n=50$.}
		\label{tab:SOCO_50}
		\begin{turn}{270}
\begin{tabular}{laaccaaccaa}
				\hline
				Problem  & \multicolumn{2}{a}{CPA-DE$_{R}^{60}$} & \multicolumn{2}{c}{DE$_{\rm bin}$} & \multicolumn{2}{a}{DE$_{\rm exp}$} & \multicolumn{2}{c}{CHC} & \multicolumn{2}{a}{GCMAES}\\
				\hline
				& Mean &  StD & Mean &  StD & Mean &  StD & Mean &  StD & Mean &  StD\\
				\hline
				$f_1$		& $5.23e-14$ & $1.57e-14$	& $3.00e-17$ & $7.69e-18$	& $2.78e-17$ & $6.29e-33$	& $2.90e+02$ & $5.69e+02$	& $2.78e-17$ & $6.29e-33$ \\
				$f_2$		& $3.33e+00$ & $1.36e+00$	& $3.87e+01$ & $8.90e+00$	& $3.31e-01$ & $5.90e-02$	& $7.72e+01$ & $1.23e+01$	& $7.69e-11$ & $4.83e-11$ \\
				$f_3$		& $4.42e+01$ & $1.77e+01$	& $6.99e+01$ & $3.58e+01$	& $3.10e+01$ & $8.65e+00$	& $5.64e+07$ & $1.42e+08$	& $6.38e-01$ & $1.49e+00$ \\
				$f_4$		& $3.98e-02$ & $1.99e-01$	& $3.21e+01$ & $1.38e+01$	& $4.79e-02$ & $2.01e-01$	& $1.12e+02$ & $2.74e+01$	& $3.72e+02$ & $8.68e+01$ \\
				$f_5$		& $2.61e-14$ & $7.87e-15$	& $9.86e-04$ & $2.76e-03$	& $0.00e+00$ & $0.00e+00$	& $9.02e-01$ & $1.82e+00$	& $2.16e-01$ & $5.64e-01$ \\
				$f_6$		& $1.18e-13$ & $1.57e-14$	& $7.16e-14$ & $1.86e-14$	& $1.39e-13$ & $9.43e-15$	& $3.23e+00$ & $2.44e+00$	& $1.90e+01$ & $1.02e+00$ \\
\hline
			\end{tabular}
		\end{turn}
	\end{center}
\end{table}
\clearpage

\begin{table}[t]
	\scriptsize
	\renewcommand{\tabcolsep}{5.0pt}
	\renewcommand{\arraystretch}{1.5}
	\begin{center}
		\caption{Average errors and standard deviations of the CPA-DE$_{R}^{60}$ and the SOCO base algorithms for dimension $n=100$.}
		\label{tab:SOCO_100}
		\begin{turn}{270}
			\begin{tabular}{laaccaaccaa}
				\hline
				Problem  & \multicolumn{2}{a}{CPA-DE$_{R}^{60}$} & \multicolumn{2}{c}{DE$_{\rm bin}$} & \multicolumn{2}{a}{DE$_{\rm exp}$} & \multicolumn{2}{c}{CHC} & \multicolumn{2}{a}{GCMAES}\\
				\hline
				& Mean &  StD & Mean &  StD & Mean &  StD & Mean &  StD & Mean &  StD\\
				\hline
				$f_1$		& $7.05e-14$ & $2.48e-14$		& $1.12e-16$ & $4.28e-17$	& $7.77e-17$ & $1.13e-17$	& $4.67e+02$ & $7.02e+02$	& $5.55e-17$ & $1.26e-32$ \\
				$f_2$		& $1.74e+01$ & $4.21e+00$		& $7.74e+01$ & $7.77e+00$	& $4.60e+00$ & $4.24e-01$	& $9.96e+01$ & $1.16e+01$	& $2.61e-03$ & $1.30e-02$ \\
				$f_3$		& $1.01e+02$ & $2.49e+01$		& $4.43e+02$ & $3.63e+02$	& $8.01e+01$ & $1.03e+01$	& $1.52e+08$ & $2.69e+08$	& $1.23e+01$ & $1.80e+01$ \\
				$f_4$		& $3.98e-02$ & $1.99e-01$		& $1.01e+02$ & $2.25e+01$	& $9.53e-03$ & $4.76e-02$	& $2.92e+02$ & $5.16e+01$	& $8.38e+02$ & $1.39e+02$ \\
				$f_5$		& $3.52e-14$ & $1.24e-14$		& $2.93e-02$ & $5.32e-02$	& $2.55e-17$ & $5.19e-18$	& $5.95e+00$ & $1.29e+01$	& $2.68e+00$ & $1.05e+01$ \\
				$f_6$		& $2.30e-13$ & $1.99e-14$		& $1.55e+00$ & $3.88e-01$	& $3.10e-13$ & $1.62e-14$	& $4.79e+00$ & $1.87e+00$	& $1.86e+01$ & $2.45e+00$ \\
\hline
			\end{tabular}
		\end{turn}
	\end{center}
\end{table}
\clearpage

\begin{table}[t]
	\scriptsize
	\renewcommand{\tabcolsep}{5.0pt}
	\renewcommand{\arraystretch}{1.5}
	\begin{center}
		\caption{Average errors and standard deviations of the CPA-DE$_{R}^{60}$ and the SOCO base algorithms for dimension $n=200$.}
		\label{tab:SOCO_200}
		\begin{turn}{270}
			\begin{tabular}{laaccaaccaa}
				\hline
				Problem  & \multicolumn{2}{a}{CPA-DE$_{R}^{60}$} & \multicolumn{2}{c}{DE$_{\rm bin}$} & \multicolumn{2}{a}{DE$_{\rm exp}$} & \multicolumn{2}{c}{CHC} & \multicolumn{2}{a}{GCMAES}\\
				\hline
				& Mean &  StD & Mean &  StD & Mean &  StD & Mean &  StD & Mean &  StD\\
				\hline
				$f_1$		& $1.59e-13$ & $2.32e-14$	& $6.39e-16$ & $5.32e-16$	& $1.78e-16$ & $1.60e-17$	& $9.61e+02$ & $1.65e+03$	& $1.17e-16$ & $1.13e-17$ \\
				$f_2$		& $3.64e+01$ & $3.41e+00$	& $1.01e+02$ & $5.90e+00$	& $1.89e+01$ & $1.05e+00$	& $1.17e+02$ & $7.60e+00$	& $7.47e-02$ & $2.44e-01$ \\
				$f_3$		& $1.96e+02$ & $2.31e+01$	& $6.38e+02$ & $3.80e+02$	& $1.79e+02$ & $8.89e+00$	& $2.54e+08$ & $3.97e+08$	& $1.24e+02$ & $8.85e+01$ \\
				$f_4$		& $3.98e-02$ & $1.99e-01$	& $4.21e+02$ & $5.63e+01$	& $8.52e-02$ & $3.98e-01$	& $6.32e+02$ & $8.43e+01$	& $1.57e+03$ & $1.54e+02$ \\
				$f_5$		& $7.28e-14$ & $1.44e-14$	& $3.00e-01$ & $7.73e-01$	& $7.49e-17$ & $6.94e-18$	& $1.02e+01$ & $1.59e+01$	& $1.13e+00$ & $2.84e+00$ \\
				$f_6$		& $4.65e-13$ & $2.58e-14$	& $5.28e+00$ & $9.65e-01$	& $6.46e-13$ & $2.53e-14$	& $8.14e+00$ & $2.26e+00$	& $1.93e+01$ & $7.39e-01$ \\
\hline
			\end{tabular}
		\end{turn}
	\end{center}
\end{table}
\clearpage

\begin{table}[t]
	\scriptsize
	\renewcommand{\tabcolsep}{5.0pt}
	\renewcommand{\arraystretch}{1.5}
	\begin{center}
		\caption{Average errors and standard deviations of the CPA-DE$_{R}^{60}$ and the SOCO base algorithms for dimension $n=500$.}
		\label{tab:SOCO_500}
		\begin{turn}{270}
			\begin{tabular}{laaccaaccaa}
				\hline
				Problem  & \multicolumn{2}{a}{CPA-DE$_{R}^{60}$} & \multicolumn{2}{c}{DE$_{\rm bin}$} & \multicolumn{2}{a}{DE$_{\rm exp}$} & \multicolumn{2}{c}{CHC} & \multicolumn{2}{a}{GCMAES}\\
				\hline
				& Mean &  StD & Mean &  StD & Mean &  StD & Mean &  StD & Mean &  StD\\
				\hline
				$f_1$		& $4.18e-13$ & $3.62e-14$	& $3.88e-05$ & $7.93e-05$	& $5.17e-16$ & $1.36e-17$	& $9.25e+02$ & $1.27e+03$	& n/a & n/a \\
				$f_2$		& $6.76e+01$ & $2.86e+00$	& $1.25e+02$ & $5.37e+00$	& $5.38e+01$ & $1.21e+00$	& $1.35e+02$ & $5.52e+00$	& n/a & n/a \\
				$f_3$		& $4.92e+02$ & $2.26e+01$	& $3.44e+04$ & $1.62e+05$	& $4.74e+02$ & $1.48e+00$	& $6.93e+08$ & $1.72e+09$	& n/a & n/a \\
				$f_4$		& $3.98e-02$ & $1.99e-01$	& $2.35e+03$ & $1.60e+02$	& $7.12e-01$ & $9.64e-01$	& $2.11e+03$ & $1.66e+02$	& n/a & n/a \\
				$f_5$		& $2.02e-13$ & $1.25e-14$	& $3.11e-01$ & $5.07e-01$	& $2.38e-16$ & $1.18e-17$	& $1.45e+01$ & $2.52e+01$	& n/a & n/a \\
				$f_6$		& $1.18e-12$ & $3.31e-14$	& $1.49e+01$ & $8.38e-01$	& $1.64e-12$ & $4.85e-14$	& $1.27e+01$ & $1.26e+00$	& n/a & n/a \\
\hline
			\end{tabular}
		\end{turn}
	\end{center}
\end{table}
\clearpage

\begin{table}[t]
	\scriptsize
	\renewcommand{\tabcolsep}{5pt}
	\renewcommand{\arraystretch}{0.85}
	\begin{center}
		\caption{Average errors between CPA-DE and a variety of algorithms on the high-dimensional SOCO 2011 test suite. Test Problems $f_{1}-f_{9}$.  (Part 1)}
		\label{tab:SOCO11_Errors1}
		\begin{turn}{270}
			\begin{tabular}{cacacacacaca}
				\hline
				& $f_1$ & $f_2$ & $f_3$ & $f_4$ & $f_5$ & $f_6$ & $f_7$ & $f_8$ & $f_9$  \\
				\hline
				CPA-DE$_{R}^{60}$\\
				$50$  	& $0.00e+00$ & $3.33e+00$ & $4.42e+01$ & $3.98e-02$ & $0.00e+00$ & $0.00e+00$ & $0.00e+00$ & $8.34e-01$ & $0.00e+00$ \\
				$100$ 	& $0.00e+00$ & $1.74e+01$ & $1.01e+02$ & $3.98e-02$ & $0.00e+00$ & $0.00e+00$ & $0.00e+00$ & $6.31e+01$ & $0.00e+00$ \\
				$200$ 	& $0.00e+00$ & $3.64e+01$ & $1.96e+02$ & $3.98e-02$ & $0.00e+00$ & $0.00e+00$ & $0.00e+00$ & $1.41e+03$ & $0.00e+00$ \\
				$500$ 	& $0.00e+00$ & $6.76e+01$ & $4.92e+02$ & $3.98e-02$ & $0.00e+00$ & $0.00e+00$ & $0.00e+00$ & $2.77e+04$ & $0.00e+00$ \\
				\hline                       
				SOUPDE	 \\                                                                                                                  
				$50$  	& $0.00e+00$ & $1.18e+00$ & $3.10e+01$ & $3.98e-02$ & $0.00e+00$ & $1.47e-14$ & $2.28e-14$ & $9.69e-02$ & $3.75e-06$ \\
\hline
			\end{tabular}
		\end{turn}
	\end{center}
\end{table}
\clearpage

\begin{table}[t]
	\scriptsize
	\renewcommand{\tabcolsep}{5pt}
	\begin{center}
		\caption{Average errors between dRL-GPALS-DE and a variety of algorithms on the high-dimensional SOCO 2011 test suite. Test Problems $f_{1}-f_{9}$. (Part 2)}
		\label{tab:SOCO11_Errors2}
		\begin{turn}{270}
			\begin{tabular}{cacacacacaca}
				\hline
				& $f_1$ & $f_2$ & $f_3$ & $f_4$ & $f_5$ & $f_6$ & $f_7$ & $f_8$ & $f_9$  \\
				\hline                       
				GODE	 \\                                                                                                                  
				$50$  	& $0.00e+00$ & $2.57e-01$ & $3.06e+01$ & $1.05e-13$ & $0.00e+00$ & $1.24e-14$ & $0.00e+00$ & $1.67e-01$ & $7.77e-06$ \\
				$100$ 	& $0.00e+00$ & $3.65e+00$ & $8.14e+01$ & $8.32e-14$ & $0.00e+00$ & $2.60e-14$ & $0.00e+00$ & $7.53e+01$ & $1.46e-05$ \\
				$200$ 	& $0.00e+00$ & $1.53e+01$ & $1.80e+02$ & $4.17e-13$ & $0.00e+00$ & $5.45e-14$ & $0.00e+00$ & $2.10e+03$ & $3.23e-05$ \\
				$500$ 	& $0.00e+00$ & $5.81e+01$ & $4.76e+02$ & $1.62e-03$ & $0.00e+00$ & $1.43e-13$ & $0.00e+00$ & $3.93e+04$ & $7.84e-05$ \\
				\hline															     
				RPSO-vm   \\                                                                                                      
				$50$  	& $2.62e-14$ & $7.54e-03$ & $1.75e+03$ & $2.30e-14$ & $9.53e-02$ & $1.49e-12$ & $1.14e-15$ & $1.06e+03$ & $2.94e-01$ \\
\hline
			\end{tabular}
		\end{turn}
	\end{center}
\end{table}
\clearpage

\begin{table}[t]
	\scriptsize
	\renewcommand{\tabcolsep}{5pt}
	\renewcommand{\arraystretch}{0.85}
	\begin{center}
		\caption{Average errors between CPA-DE and a variety of algorithms on the high-dimensional SOCO 2011 test suite. Test Problems $f_{10}-f_{19}$.  (Part 1)}
		\label{tab:SOCO11_Errors3}
		\begin{turn}{270}
			\begin{tabular}{cacacacacaca}
				\hline
				& $f_{10}$ &  $f_{11}$ & $f_{12}$ & $f_{13}$ & $f_{14}$ & $f_{15}$ & $f_{16}$ & $f_{17}$ & $f_{18}$ & $f_{19}$ \\
				\hline
				CPA-DE$_{R}^{60}$\\
				$50$   & $0.00e+00$ & $0.00e+00$ & $0.00e+00$ & $2.74e+01$ & $0.00e+00$ & $0.00e+00$ & $0.00e+00$ & $1.22e+00$ & $1.46e-08$ & $0.00e+00$ \\
				$100$  & $0.00e+00$ & $0.00e+00$ & $0.00e+00$ & $7.26e+01$ & $3.98e-02$ & $0.00e+00$ & $0.00e+00$ & $2.43e+01$ & $3.56e-08$ & $0.00e+00$ \\
				$200$  & $0.00e+00$ & $0.00e+00$ & $0.00e+00$ & $1.51e+02$ & $0.00e+00$ & $0.00e+00$ & $0.00e+00$ & $5.19e+01$ & $3.98e-02$ & $0.00e+00$ \\
				$500$  & $0.00e+00$ & $0.00e+00$ & $0.00e+00$ & $3.76e+02$ & $0.00e+00$ & $0.00e+00$ & $0.00e+00$ & $1.26e+02$ & $2.77e-07$ & $0.00e+00$ \\
				\hline
				SOUPDE \\
				$50$   & $0.00e+00$ & $3.09e-06$ & $0.00e+00$ & $2.06e+01$ & $0.00e+00$ & $1.38e-14$ & $0.00e+00$ & $2.53e-01$ & $0.00e+00$ & $0.00e+00$ \\
\hline
			\end{tabular}
		\end{turn}
	\end{center}
\end{table}
\clearpage

\begin{table}[t]
	\scriptsize
	\renewcommand{\tabcolsep}{5pt}
	\begin{center}
		\caption{Average errors between dRL-GPALS-DE and a variety of algorithms on the high-dimensional SOCO 2011 test suite. Test Problems $f_{10}-f_{19}$. (Part 2)}
		\label{tab:SOCO11_Errors4}
		\begin{turn}{270}
			\begin{tabular}{cacacacacaca}
				\hline
				& $f_{10}$ &  $f_{11}$ & $f_{12}$ & $f_{13}$ & $f_{14}$ & $f_{15}$ & $f_{16}$ & $f_{17}$ & $f_{18}$ & $f_{19}$ \\
				\hline
				GODE \\
				$50$   & $0.00e+00$ & $6.44e-06$ & $1.33e-13$ & $2.55e+01$ & $6.24e-09$ & $0.00e+00$ & $1.57e-10$ & $1.17e+00$ & $2.97e-07$ & $0.00e+00$ \\
				$100$  & $0.00e+00$ & $1.58e-05$ & $7.57e-12$ & $6.32e+01$ & $4.13e-08$ & $0.00e+00$ & $3.75e-10$ & $1.11e+01$ & $1.11e-06$ & $0.00e+00$ \\
				$200$  & $0.00e+00$ & $3.12e-05$ & $1.20e-10$ & $1.38e+02$ & $8.17e-02$ & $0.00e+00$ & $9.54e-10$ & $3.74e+01$ & $1.91e-06$ & $0.00e+00$ \\
				$500$  & $0.00e+00$ & $8.25e-05$ & $7.39e-10$ & $3.59e+02$ & $7.67e-02$ & $0.00e+00$ & $2.24e-09$ & $1.12e+02$ & $5.06e-06$ & $0.00e+00$ \\
				\hline                                      
				RPSO-vm \\                                                                                                                              
				$50$   & $0.00e+00$ & $1.68e-02$ & $8.58e-02$ & $6.57e+02$ & $6.81e-02$ & $0.00e+00$ & $7.88e-11$ & $8.73e+02$ & $5.05e-02$ & $0.00e+00$ \\
\hline
			\end{tabular}
		\end{turn}
	\end{center}
\end{table}
\clearpage

\begin{table}[t]
	\scriptsize
	\renewcommand{\tabcolsep}{5pt}
	\renewcommand{\arraystretch}{0.85}
	\begin{center}
		\caption{Average errors of CPA-DE and SHADE variants for the high-dimensional SOCO 2011 test suite. Test problems $f_{1}-f_{9}$.}
		\label{tab:SHADE-SOCO-ERRORS1}
		\begin{turn}{270}
			\begin{tabular}{cacacacacaca}
				\hline
				& $f_1$ & $f_2$ & $f_3$ & $f_4$ & $f_5$ & $f_6$ & $f_7$ & $f_8$ & $f_9$  \\
				\hline 
				CPA-DE$_{R}^{60}$ \\                                                                                                      
				$50$  	& $0.00e+00$ & $3.33e+00$ & $4.42e+01$ & $3.98e-02$ & $0.00e+00$ & $0.00e+00$ & $0.00e+00$ & $8.34e-01$ & $0.00e+00$ \\
				$100$ 	& $0.00e+00$ & $1.74e+01$ & $1.01e+02$ & $3.98e-02$ & $0.00e+00$ & $0.00e+00$ & $0.00e+00$ & $6.31e+01$ & $0.00e+00$ \\
				$200$ 	& $0.00e+00$ & $3.64e+01$ & $1.96e+02$ & $3.98e-02$ & $0.00e+00$ & $0.00e+00$ & $0.00e+00$ & $1.41e+03$ & $0.00e+00$ \\
				$500$ 	& $0.00e+00$ & $6.76e+01$ & $4.92e+02$ & $3.98e-02$ & $0.00e+00$ & $0.00e+00$ & $0.00e+00$ & $2.77e+04$ & $0.00e+00$ \\
				\hline  
				
				CPA-DE$_{C}^{60}$ \\                                                                                                      
				$50$  	& $0.00e+00$ & $1.57e+00$ & $4.77e+00$ & $5.23e+01$ & $5.80e-03$ & $5.07e-01$ & $0.00e+00$ & $8.50e-03$ & $1.97e+00$ \\
\hline 															     
			\end{tabular}
		\end{turn}
	\end{center}
\end{table}
\clearpage

\begin{table}[t]
	\scriptsize
	\renewcommand{\tabcolsep}{5pt}
	\renewcommand{\arraystretch}{0.85}
	\begin{center}
		\caption{Average errors of CPA-DE and SHADE variants for the high-dimensional SOCO 2011 test suite. Test problems $f_{10}-f_{19}$.}
		\label{tab:SHADE-SOCO-ERRORS2}
		\begin{turn}{270}
			\begin{tabular}{cacacacacaca}
				\hline
				& $f_{10}$ &  $f_{11}$ & $f_{12}$ & $f_{13}$ & $f_{14}$ & $f_{15}$ & $f_{16}$ & $f_{17}$ & $f_{18}$ & $f_{19}$ \\
				\hline
				CPA-DE$_{R}^{60}$ \\                                                                                                        
				$50$   & $0.00e+00$ & $0.00e+00$ & $0.00e+00$ & $2.74e+01$ & $0.00e+00$ & $0.00e+00$ & $0.00e+00$ & $1.22e+00$ & $1.46e-08$ & $0.00e+00$ \\
				$100$  & $0.00e+00$ & $0.00e+00$ & $0.00e+00$ & $7.26e+01$ & $3.98e-02$ & $0.00e+00$ & $0.00e+00$ & $2.43e+01$ & $3.56e-08$ & $0.00e+00$ \\
				$200$  & $0.00e+00$ & $0.00e+00$ & $0.00e+00$ & $1.51e+02$ & $0.00e+00$ & $0.00e+00$ & $0.00e+00$ & $5.19e+01$ & $3.98e-02$ & $0.00e+00$ \\
				$500$  & $0.00e+00$ & $0.00e+00$ & $0.00e+00$ & $3.76e+02$ & $0.00e+00$ & $0.00e+00$ & $0.00e+00$ & $1.26e+02$ & $2.77e-07$ & $0.00e+00$ \\
				\hline  
				
				CPA-DE$_{C}^{60}$ \\                                                                                                                       
				$50$   & $8.07e+00$ & $2.43e+00$ & $2.86e+01$ & $7.60e+01$ & $3.92e+01$ & $1.61e-01$ & $2.68e+01$ & $3.31e+01$ & $1.86e-01$ & $2.99e+00$ \\
\hline
			\end{tabular}
		\end{turn}
	\end{center}
\end{table}
\clearpage


\section{Results for CEC 2013}
\label{sec:app2}

\begin{table}[t]
	\scriptsize
	\renewcommand{\tabcolsep}{5pt}
	\renewcommand{\arraystretch}{0.6}
	\begin{center}
		\caption{Average errors between CPA-DE and a variety of algorithms on the low-dimensional CEC 2013 test suite. Test Problems $f_{1}-f_{9}$.}
		\label{tab:CEC13-ERRORS1}
		\begin{turn}{270}
\begin{tabular}{cacacacacac}
				\hline
				& $f_1$ & $f_2$ & $f_3$ & $f_4$ & $f_5$ & $f_6$ & $f_7$ & $f_8$ & $f_9$  \\
				\hline
				CPA-DE$_{C}^{60}$\\
				$30$  	& $0.00e+00$ & $7.90e+05$ & $2.50e+06$ & $3.38e+03$ & $0.00e+00$ & $1.93e+01$ & $3.70e+01$ & $2.09e+01$ & $1.70e+01$ \\
				$50$ 	& $0.00e+00$ & $2.25e+06$ & $3.72e+07$ & $9.52e+03$ & $0.00e+00$ & $4.45e+01$ & $8.09e+01$ & $2.11e+01$ & $4.02e+01$ \\
				\hline
				
				NIPOPaCMA\\
				$30$  	& $0.00e+00$ & $0.00e+00$ & $7.45e-07$ & $0.00e+00$ & $0.00e+00$ & $0.00e+00$ & $4.05e+00$ & $2.09e+01$ & $2.82e+00$ \\
				$50$ 	& $0.00e+00$ & $0.00e+00$ & $1.88e+01$ & $0.00e+00$ & $0.00e+00$ & $0.00e+00$ & $1.17e+01$ & $2.11e+01$ & $6.88e+00$ \\
				\hline

				icmaesils\\
\hline
			\end{tabular}
		\end{turn}
	\end{center}
\end{table}
\clearpage

\begin{table}[t]
	\scriptsize
	\renewcommand{\tabcolsep}{5pt}
	\renewcommand{\arraystretch}{0.6}
	\begin{center}
		\caption{Average errors between CPA-DE and a variety of algorithms on the low-dimensional CEC 2013 test suite. Test Problems $f_{10}-f_{19}$.}
		\label{tab:CEC13-ERRORS2}
		\begin{turn}{270}
			\begin{tabular}{cacacacacac}
				\hline
				& $f_{10}$ &  $f_{11}$ & $f_{12}$ & $f_{13}$ & $f_{14}$ & $f_{15}$ & $f_{16}$ & $f_{17}$ & $f_{18}$ & $f_{19}$ \\
				\hline
				CPA-DE$_{C}^{60}$ \\
				$30$   & $8.91e-02$ & $0.00e+00$ & $6.64e+01$ & $1.02e+02$ & $3.06e-02$ & $3.16e+03$ & $1.76e+00$ & $3.04e+01$ & $9.97e+01$ & $1.11e+00$  \\
				$50$   & $1.83e-01$ & $0.00e+00$ & $1.63e+02$ & $2.31e+02$ & $9.77e-02$ & $6.77e+03$ & $2.25e+00$ & $5.08e+01$ & $2.05e+02$ & $2.07e+00$ \\
				\hline
				
				NIPOPaCMA\\
				$30$   & $0.00e+00$ & $1.03e+00$ & $6.57e-01$ & $9.32e-01$ & $7.09e+02$ & $6.70e+02$ & $2.48e+00$ & $3.42e+01$ & $5.40e+01$ & $2.41e+00$ \\
				$50$   & $0.00e+00$ & $1.99e+00$ & $1.37e+00$ & $1.48e+00$ & $1.27e+03$ & $1.35e+03$ & $3.37e+00$ & $5.77e+01$ & $1.94e+02$ & $4.47e+00$ \\
				\hline

				icmaesils\\
\hline
			\end{tabular}
		\end{turn}
	\end{center}
\end{table}
\clearpage

\begin{table}[t]
	\scriptsize
	\renewcommand{\tabcolsep}{5pt}
	\renewcommand{\arraystretch}{0.6}
	\begin{center}
		\caption{Average errors between CPA-DE and a variety of algorithms on the low-dimensional CEC 2013 test suite. Test Problems $f_{20}-f_{28}$.}
		\label{tab:CEC13-ERRORS3}
		\begin{turn}{270}
			\begin{tabular}{cacacacaca}
				\hline
				& $f_{20}$ & $f_{21}$ & $f_{22}$ & $f_{23}$ & $f_{24}$ & $f_{25}$ & $f_{26}$ & $f_{27}$ & $f_{28}$  \\
				\hline
				CPA-DE$_{C}^{60}$\\
				$30$  	& $1.15e+01$ & $3.23e+02$ & $7.38e+01$ & $3.82e+03$ & $2.39e+02$ & $2.65e+02$ & $2.00e+02$ & $7.65e+02$ & $3.00e+02$ \\
				$50$ 	& $2.13e+01$ & $5.15e+02$ & $1.75e+01$ & $8.16e+03$ & $2.98e+02$ & $3.37e+02$ & $2.00e+02$ & $1.36e+03$ & $4.00e+02$ \\
				\hline
				
				NIPOPaCMA\\
				$30$  	& $1.34e+01$ & $2.41e+02$ & $5.75e+02$ & $6.67e+02$ & $2.91e+02$ & $2.79e+02$ & $2.51e+02$ & $8.70e+02$ & $3.00e+02$ \\
				$50$ 	& $2.30e+01$ & $3.65e+02$ & $1.02e+03$ & $1.19e+03$ & $3.70e+02$ & $3.65e+02$ & $2.88e+02$ & $1.90e+03$ & $5.72e+02$ \\
				\hline   

				icmaesils\\
\hline
			\end{tabular}
		\end{turn}
	\end{center}
\end{table}
\clearpage

\begin{table}[t]
	\scriptsize
	\renewcommand{\tabcolsep}{5pt}
	\begin{center}
		\caption{Average errors of CPA-DE and SHADE variants for the low-dimensional CEC 2013 test suite. Test Problems $f_{1}-f_{9}$.}
		\label{tab:SHADE-CEC13-ERRORS1}
		\begin{turn}{270}
			\begin{tabular}{cacacacacac}
				\hline
				& $f_1$ & $f_2$ & $f_3$ & $f_4$ & $f_5$ & $f_6$ & $f_7$ & $f_8$ & $f_9$  \\
				\hline
				CPA-DE$_{R}^{60}$ \\
				$30$  	& $0.00e+00$ & $7.90e+05$ & $2.50e+06$ & $3.38e+03$ & $0.00e+00$ & $1.93e+01$ & $3.70e+01$ & $2.09e+01$ & $1.70e+01$ \\
				$50$ 	& $0.00e+00$ & $2.25e+06$ & $3.72e+07$ & $9.52e+03$ & $0.00e+00$ & $4.45e+01$ & $8.09e+01$ & $2.11e+01$ & $4.02e+01$ \\
				\hline
				
				CPA-DE$_{C}^{60}$ \\   
				$30$  	& $0.00e+00$ & $7.74e+04$ & $5.31e+06$ & $2.92e-01$ & $0.00e+00$ & $4.00e+00$ & $1.50e+01$ & $2.09e+01$ & $2.39e+01$ \\
				$50$ 	& $0.00e+00$ & $5.45e+05$ & $4.12e+07$ & $3.86e+00$ & $0.00e+00$ & $4.17e+01$ & $5.32e+01$ & $2.11e+01$ & $4.69e+01$ \\
				\hline                          
				
				CPA-DE$_{R}^{100}$ \\
\hline			
				
			\end{tabular}
		\end{turn}
	\end{center}
\end{table}
\clearpage

\begin{table}[t]
	\scriptsize
	\renewcommand{\tabcolsep}{5pt}
	\begin{center}
		\caption{Average errors of CPA-DE and SHADE variants for the low-dimensional CEC 2013 test suite. Test Problems $f_{10}-f_{19}$.}
		\label{tab:SHADE-CEC13-ERRORS2}
		\begin{turn}{270}
			\begin{tabular}{cacacacacac}
				\hline
				& $f_{10}$ &  $f_{11}$ & $f_{12}$ & $f_{13}$ & $f_{14}$ & $f_{15}$ & $f_{16}$ & $f_{17}$ & $f_{18}$ & $f_{19}$ \\
				\hline
				CPA-DE$_{R}^{60}$ \\
				$30$  & $8.91e-02$ & $0.00e+00$ & $6.64e+01$ & $1.02e+02$ & $3.06e-02$ & $3.16e+03$ & $1.76e+00$ & $3.04e+01$ & $9.97e+01$ & $1.11e+00$ \\
				$50$  & $1.83e-01$ & $0.00e+00$ & $1.63e+02$ & $2.31e+02$ & $9.77e-02$ & $6.77e+03$ & $2.25e+00$ & $5.08e+01$ & $2.05e+02$ & $2.07e+00$ \\
				
				\hline
				CPA-DE$_{C}^{60}$ \\  
				$30$  & $5.40e-02$ & $2.83e-01$ & $4.20e+01$ & $7.91e+01$ & $8.06e+02$ & $4.62e+03$ & $1.75e+00$ & $4.58e+01$ & $1.26e+02$ & $3.24e+00$ \\
				$50$  & $7.71e-02$ & $1.41e+01$ & $1.09e+02$ & $2.23e+02$ & $2.28e+03$ & $9.85e+03$ & $2.30e+00$ & $9.97e+01$ & $2.72e+02$ & $1.04e+01$ \\
				\hline
				
				CPA-DE$_{R}^{100}$ \\
\hline
			\end{tabular}
		\end{turn}
	\end{center}
\end{table}
\clearpage

\begin{table}[t]
	\scriptsize
	\renewcommand{\tabcolsep}{5pt}
	\begin{center}
		\caption{Average errors of CPA-DE and SHADE variants for the low-dimensional CEC 2013 test suite. Test Problems $f_{20}-f_{28}$.}
		\label{tab:SHADE-CEC13-ERRORS3}
		\begin{turn}{270}
			\begin{tabular}{cacacacaca}
				\hline
				& $f_{20}$ & $f_{21}$ & $f_{22}$ & $f_{23}$ & $f_{24}$ & $f_{25}$ & $f_{26}$ & $f_{27}$ & $f_{28}$  \\
				\hline
				CPA-DE$_{R}^{60}$ \\
				$30$  	& $1.15e+01$ & $3.23e+02$ & $7.38e+01$ & $3.82e+03$ & $2.39e+02$ & $2.65e+02$ & $2.00e+02$ & $7.65e+02$ & $3.00e+02$ \\
				$50$ 	& $2.13e+01$ & $5.15e+02$ & $1.75e+01$ & $8.16e+03$ & $2.98e+02$ & $3.37e+02$ & $2.00e+02$ & $1.36e+03$ & $4.00e+02$ \\
				\hline
				
				CPA-DE$_{C}^{60}$ \\ 
				$30$  	& $1.11e+01$ & $3.28e+02$ & $7.49e+02$ & $4.82e+03$ & $2.23e+02$ & $2.54e+02$ & $2.33e+02$ & $5.44e+02$ & $3.41e+02$ \\
				$50$ 	& $2.04e+01$ & $8.45e+02$ & $2.06e+03$ & $9.95e+03$ & $2.68e+02$ & $3.13e+02$ & $3.11e+02$ & $1.00e+03$ & $5.80e+02$ \\
				\hline                          
				
				CPA-DE$_{R}^{100}$ \\
\hline
			\end{tabular}
		\end{turn}
	\end{center}
\end{table}
\clearpage


\end{document}